\documentclass[a4paper,10pt]{article}
\usepackage{amsmath}
\usepackage{amssymb}
\usepackage{bbm}
\usepackage[pdftex]{graphicx}
\usepackage{color}

\numberwithin{equation}{section}

\newtheorem{thm}{Theorem}[section]
\newtheorem{cor}[thm]{Corollary}
\newtheorem{lem}[thm]{Lemma}
\newtheorem{prop}[thm]{Proposition}

\newtheorem{defn}{Definition}[section]
\newtheorem{rem}{Remark}[section]

\newcommand{\ii}[2][\delta']{\left\langle #2 \right
\rangle^{#1}}
\def\conv{\mathrm{conv}}
\def\GM{G_{\max}}
\def\Gm{G_{\min}}

\def\d{\mathtt{d}}

\def\ll{\left\lbrace}
\def\rr{\right\rbrace}

\def\lambdam{\lambda_{\min}}
\def\lambdaM{\lambda_{\max}}
\def\kappam{\kappa_{\min}}
\def\kappaM{\kappa_{\max}}

\title{Duality of Floating and Illumination Bodies
}

\author{Olaf Mordhorst \thanks{Partially supported by the German Academic Exchange Service} and  Elisabeth M.  Werner \thanks{Partially supported by  NSF grant DMS-1504701}}

\date{}
\begin{document}

\maketitle

\begin{abstract}
We investigate a duality relation between floating and illumination bodies. The definitions of these two bodies suggest that the polar of the floating body should be similar to the illumination body of the polar. We consider this question for the class of centrally symmetric convex bodies. We provide precise estimates for \(B_p^n\) and for centrally symmetric convex bodies with everywhere positive Gauss curvature. Our estimates show that equality of the polar of the floating body and the illumination body of the polar can only be achieved in the case of ellipsoids.
\end{abstract}

\section{Introduction}
Floating bodies  and illumination bodies  are attracting considerable interest as 
their important properties  make them effective and powerful tools. Therefore they,  and  
  the related affine surface areas,  are
omnipresent in geometry, e.g.,  \cite{HaberlSchuster:2009, Haberl:2014, LR:2010, Lutwak:1991, Lutwak:1996,  HuangLutwakYangZhang, MW:2000, Schuett:2004, WernerYe11}
and  find applications in  many other areas such as
information theory, e.g., \cite{Artstein-Avidan:2012, PaourisWerner2012, Werner2012}, 
the study  of polytopes and approximation by polytopes 
\cite{Barany:1992, Boeroetzky2000, Boeroetzky2000/2, GruberBook,  Ludwig:1998, Reitzner:2005, Reitzner:2010, Schuett91,Schuett:1994,  Schuett:2003}
and 
partial differential equations (e.g., \cite{Lutwak:1995, TW4} and  the solutions  for the affine Bernstein and  Plateau problems
by Trudinger and Wang \cite {TW1,TW2, TW3}). 
\par
Very recent developments are the introduction of  the floating body in spherical space \cite{BW:2015}
and in hyperbolic space  \cite{BW:2016}. This has already given rise to applications in approximation of spherical and hyperbolic 
convex bodies by polytopes  \cite{BLW:2016} . 
\par
A notion of floating body appeared already in the work of C.\ Dupin \cite{Dupin:1822} in 1822.
In 1990, a new definition was given by Sch\"utt and Werner \cite{SchuettWerner1990} and independently by B{\'a}r{\'a}ny and Larman \cite{Barany:1988}.
They introduced the \emph{convex} floating body as the intersection of all halfspaces whose hyperplanes
cut off a set of fixed volume of a convex body (a compact convex set).
In contrast to the original definition, the convex floating body is always convex and coincides with Dupin's floating
body if it exists.
\par
The illumination body was introduced in \cite{Werner94} as the set of those points whose convex hull with a given convex body  have fixed volume.
\vskip 2mm
The definitions of the floating body and the illumination body suggest a possible duality relation, namely that the polar of a floating body 
of a convex body $K$  is ``close" to  an illumination body of the polar of  $K$.  In fact,  for the Euclidean unit ball \(B_2^n\),  equality can always be achieved.
Note  however that equality cannot be achieved in  general  since 
it was shown in \cite{SchuettWerner1990} that floating bodies are always strictly convex,  but the illumination body of a polytope is always a polytope. 
\par
In this paper we clarify the  duality relation between floating body and illumination in the case of centrally symmetric convex bodies with \(C^2\)-boundary. We provide asymptotically precise formulas for bodies with everywhere strictly positive Gauss curvature and for \(B_p^n\), \(2\leq p<\infty\). 
\par
Floating bodies and illumination  bodies 
allow to establish  the long sought extensions of an  important affine invariant, the  affine surface area,
to general convex bodies in all
dimensions. This was carried out in \cite{SchuettWerner1990}, respectively \cite{Werner94}. 
In both instances, affine surface area appears as a limit of the volume difference of the convex body and its floating body, respectively illumination body.
Other extensions - all  coincide - were given by Lutwak \cite{Lutwak:1996} and Leichtweiss \cite{Leichtweiss:1986}.
\par
Here we carry out a limit procedure which  leads to a new affine invariant that is  different from the affine surface area. 
It is related  to the cone measure of the convex body.  These measures play a central role in many aspects of convex geometry,  e.g.,  \cite{BLYZ, CLM, Naor, PaourisWerner2012}.
\vskip 2mm
In a forthcoming paper we investigate the case of centrally symmetric polytopes which leads to discrete versions of the formulas we derive in the following.
\vskip 2mm
The  paper is organized as follows. In the next subsection we introduce  notation that is used throughout the paper and present the main theorems and some  consequences, a characterization of ellipsoids among them. Section \ref{SectionBackground} provides the necessary definitions and 
background. In Section \ref{UpperBound} we give an upper bound for general centrally symmetric convex bodies that are \(C^2\) and  in Section \ref{SectionLowerBound} we prove that this upper bound is precise in the case of \(C_+^2\)-boundary and for \(B_p^n\), \(2\leq p<\infty\). We also provide a lower bound for the case \(1<p<2\). 

\subsection{Notation}
We denote by   \(\mathbb{R}_{\geq 0}\) the non-negative real numbers. 
A convex body   \(K\subseteq \)  in \( \mathbb{R}^n\) is a  convex, compact subset of \( \mathbb{R}^n\) with non-empty interior.  \(K\) is called centrally symmetric (with respect to the origin) if \(K=-K\). From now on,  we will always denote by \(C\subseteq \mathbb{R}^n\) a centrally symmetric convex body and by   \(S\subseteq \mathbb{R}^n\) a centrally symmetric convex body  with \(C^2\)-boundary.
We refer to the books by Schneider \cite{SchneiderBook} or Gardner \cite{GardnerBook} for background on convex bodies.
\par
Let $A$ and $B$ be subsets of $\mathbb{R}^n$. Then $\mathrm{conv}[A,B] = \{ \lambda a + (1-\lambda) b: a \in A, b \in B, 0 \leq \lambda \leq 1\}$ is the convex hull of $A$ and $B$. If $B=\{x\}$, we simply write $\mathrm{conv}[A,x]$.
For a measurable set $A \in \mathbb{R}^n$, we denote by $|A|_k$ its $k$-dimensional Hausdorff measure and,  in particular, by \(|A|_n\) its \(n\)-dimensional volume.
\par
For $1\leq p < \infty$, let  $B^n_p$ be the unit ball of the space \( \mathbb{R}^n\)  equipped with the norm $\|x\|_p = \left(\sum _{i=1}^n|x_i|^p \right)^\frac{1}{p}$, 
$$B^n_p=\{x  \in\mathbb{R}^n:
\|x\|_p  \leq 1\}. $$ 

\subsection{Statement of principal results}

Let \(K\) be  convex body  in \( \mathbb{R}^n\) and \(\delta \geq 0\). The convex floating body $K_\delta$ of $K$
was introduced in \cite{SchuettWerner1990} and independently by B{\'a}r{\'a}ny and Larman \cite{Barany:1988} as the intersection of all half spaces whose defining hyperplanes cut off a set of volume $\delta |K|_n$ from $K$.
More precisely, 
\begin{equation}\label{floatingbody}
K_{\delta}= \bigcap\limits_{|K \cap H^{-}|_n\leq \delta |K|_n} H^+ , 
\end{equation}
where \(H\) is a hyperplane and \(H^+, H^-\) are the corresponding closed halfspaces. 
An important result by Meyer and Reisner  \cite{MeyerReisner} which we will use throughout, states that for centrally symmetric convex bodies 
Dupin's floating body always exists and coincides with the convex floating body. 
\vskip 2mm
The  illumination body $K^\delta$ of  $K$  was introduced in \cite{Werner94} as follows  
\begin{equation}\label{illuminationbody}
K^{\delta}=\{x\in\mathbb{R}^n: |\mathrm{conv}[K,x]|_n\leq (1+\delta)|K|_n\}\quad.
\end{equation}
Note that the illumination body is always convex. This can easily be seen by the fact that
\[
|\mathrm{conv}[K,x]|_n=\frac{1}{2}\left(|K|_n+\frac{1}{n}\int_{\partial K}|\langle x-y, u(y)\rangle|\mathrm{d}\mu(y)\right), 
\]
where $\langle \cdot,  \cdot \rangle$ is the standard inner product on $\mathbb{R}^n$, \(\mu\) is the surface measure on $\partial K$, the boundary of $K$,  and \(u(y)\) the almost everywhere uniquely determined outer normal at \(y\in \partial K\).
\vskip 2mm
The definitions of the floating body and the illumination body suggest a  duality relation, namely that the polar of a  floating body of a convex body $K$  is ``close" to the illumination body of the polar of $K$, 
\begin{equation}\label{gleich}
\left[K_{\delta}\right]^{\circ}\approx\left[K^{\circ}\right]^{\delta'}, 
\end{equation}
for suitable \(\delta\) and \(\delta'\). Note  that equality cannot be achieved in  general since 
it was shown in \cite{SchuettWerner1990} that floating bodies are always strictly convex,  but the illumination body of a polytope is again a polytope. On the other hand, equality can always be achieved for the Euclidean unit ball \(B_2^n\).
We make  (\ref{gleich}) precise in the case of  centrally symmetric convex bodies. To do so, we use the distance which  we introduce next.
\vskip 2mm
For a convex body \(K\subseteq \mathbb{R}^n\) with \(0\in K\) and \(x\in\mathbb{R}^n\backslash\{0\}\) we denote by \(r_K(x)=\max\{\lambda\geq 0: \lambda x\in K\}\) the radial function of \(K\).  We define a distance \(\d\) on the set of \(n\)-dimensional convex bodies which compares radial functions.  We only consider this distance for centrally symmetric convex bodies. If  \(C_1\) and \(C_2\)  are  \(n\)-dimensional  centrally symmetric convex bodies, we define
We define
\begin{equation}\label{distance}
\d(C_1,C_2)=\sup_{u\in\mathbb{S}^{n-1}}\max\left[\frac{r_{C_1}(u)}{r_{C_2}(u)},\frac{r_{C_2}(u)}{r_{C_1}(u)}\right]=\inf \left\{a\geq 1: \frac{1}{a}C_1\subseteq C_2 \subseteq a C_1\right\}\quad.
\end{equation}
Note that \(\log \d(\cdot,\cdot)\) is a metric which induces the same topology as the Hausdorff distance. We put 
\[
\ii[\delta]{C}= \left(\left[C^{\circ}\right]^{\delta}\right)^{\circ}\quad.
\]
\begin{defn} \label{FloatIlluMeas}
Let \(C\) be a centrally symmetric convex body and \(0<\delta<\frac{1}{2}\). We define 
\[
\d_C(\delta)=\inf_{\delta'>0}\d\left(C_{\delta},\ii{C}\right)\quad.
\] 
If $C=B_p^n$, we put \(\d_p(\delta)=\d_{B_p^n}(\delta)\).
\end{defn}
Please note that \(\d_{L(C)}(\delta)=\d_C(\delta)\) for every linear invertible map \(L\) and that $\d_2(\delta)=\d_{B_2^n}(\delta) = 1$.
\vskip 3mm
Together with the distance  the distance $ \d$, the following expressions  will be crucial to make  the relation (\ref{gleich}) precise. 
For a centrally symmetric convex body \(C\) and \(x\in\partial C\) with a unique outer normal \(u(x)\) and such that the Gauss curvature \(\kappa(x)\) exists,  we set
\begin{eqnarray} \label{G}
G_C(x)=c(C,n)\frac{\kappa(x)^{\frac{1}{n+1}}}{\langle x, u(x)\rangle}\quad, 
\end{eqnarray}
where \(c(C,n)=\frac{(n+1)^{\frac{2}{n+1}}}{2}\left(\frac{|C|_n}{|B_2^{n-1}|_{n-1}}\right)^{\frac{2}{n+1}}\). In most cases we omit the subscript \(C\) in \(G_C\). 
\par
Our two main theorems treat the case that the centrally symmetric convex body  has $C^2$-boundary.  If in addition the Gauss curvature is strictly positive everywhere,  we say that a convex body is of class \(C_+^2\). For such bodies the function \(G:\partial S\rightarrow\mathbb{R}_{\geq 0}\) is continuous with respect to the Euclidean distance. We put
$$
G_{\max} =\max_{x\in\partial S}G(x)  \hskip 3mm \text{   and  } \hskip 3mm    G_{\min}=\min_{x\in\partial S}G(x).
$$
\begin{thm} \label{main1}
\vskip 2mm
\noindent
Let $S \subset \mathbb{R}^n$ be a centrally symmetric convex body that is of class  $C^2$. Then
$$
 \limsup _{\delta \rightarrow 0} \frac{\d_S(\delta) -1}{\delta^{\frac{2}{n+1}} }
 \leq  \GM -\Gm\quad.
$$
\end{thm}
\begin{thm} \label{main2}
\vskip 2mm
\noindent
Let $S \subset \mathbb{R}^n$ be a centrally symmetric convex body that  is of class  $C^2$. If \(S\) has everywhere strictly positive Gauss curvature or is an \(B_p^n\)-ball, \(2\leq p<\infty\), then
$$
 \lim _{\delta \rightarrow 0} \frac{\d_S(\delta) -1}{\delta^{\frac{2}{n+1}} }
 =  \GM -\Gm\quad.
$$
\end{thm}
One might ask if Theorem \ref{main2} holds for general \(S\) with \(C^2\)-boundary. However, the authors think that this is not the case and that  they can construct a counterexample.
\vskip 2mm
\noindent
To put the above theorems into context, we recall that it 
was shown in \cite{SchuettWerner1990}  that for all convex bodies $K$, 
\begin{align}
\frac{1}{c_n}    \ \lim_{\delta\rightarrow 0}
\frac{|K|_n-|K_{\delta/|K|_n}|_n}{\delta^{\frac{2}{n+1}}}=  \int_{\partial K}\kappa^{\frac{1}{n+1}}(x)\mathrm{d}\mu(x)\quad,\label{SurfaceAreaFormula}
\end{align}
where $c_n= \frac{1}{2}\left(\frac{n+1 }{|B_2^{n-1}|_{n-1}}\right)^{\frac{2}{n+1}}$ and  the right-hand side integral  is the well-known affine surface area.  Note that a different normalization is chosen for the parameter \(\delta\)  in the definition of the floating body in \cite{SchuettWerner1990}. 
\par
For $x \in \partial K$ with outer normal $u(x)$, let 
 $m_{K} (x) = \frac{1}{n}  \langle x, u(x)\rangle $ be  the density function of the cone measure $M_K$ of $K$  with respect to the surface measure of \(K\). For a Borel set $A \in \partial K$,  it is  defined as $M_K(A) = | \conv[0, A]|_n$.  We write  $n_{K} (x) = \frac{1}{n |K|_n}  \langle x, u(x)\rangle $ for  the density of the normalized cone measure $\mathbb{P}_K$ of $K$. 
This means  that 
$$
d M_K (x) = m_{K} (x) d   \mu_K(x) \  \   \text {and } \  \    d \mathbb{P}_K (x) = n_{K} (x) d   \mu_K(x).
$$
If we rewrite  (\ref{SurfaceAreaFormula}) using the cone measure equivalently in such a way  that have  both sides are affine invariants, we get 
\begin{align}
 \lim_{\delta\rightarrow 0}\frac{|K|_n-|K_{\delta}|_n}{n|K|_n\delta^{\frac{2}{n+1}}}=\int_{\partial K}G(x)\mathrm{d}n_K(x)=\|G\|_{1,n_s}\label{SurfaceAreaFormulaII} . 
\end{align}
One can  consider other differences than the  volume difference on the left-hand side of (\ref{SurfaceAreaFormulaII}). Indeed, combining Proposition \ref{LowerBoundFloatingBodyC2Case} and Lemma \ref{FloatingBodyCurvature} of the following sections we get a radial-version of (\ref{SurfaceAreaFormulaII}) for
centrally symmetric  \(C^2\)-bodies
\[
\lim_{\delta\rightarrow 0}\sup_{u\in\mathbb{S}^{n-1}}\frac{r_{S}(u)-r_{S_{\delta}}(u)}{r_S(u)\delta^{\frac{2}{n+1}}}=\lim_{\delta\rightarrow 0}\frac{\d(S_{\delta},S)}{\delta^{\frac{2}{n+1}}}=\GM=\|G\|_{\infty, n_s} \quad.
\]
\vskip 2mm
Let $S$ be a centrally symmetric convex body that is $C^2_+$. There is  a nice way to write \(G_S\)  of (\ref{G}) in terms of cone measures. 
Let  $m_S$ and $n_s$ be the cone measure respectively the normalized cone measure of $S$  defined above and let $c_n$  be as above. Observe that 
\begin{eqnarray} \label{conemeasure}
G_S(x)&=&c_n \   |S|_n  ^{\frac{2}{n+1}} \frac{\kappa(x)^{\frac{1}{n+1}}}{\langle x, u(x)\rangle} = c_n \  |S|_n ^{\frac{2}{n+1}} \left(\frac{\kappa(x)}{\langle x, u(x)\rangle^n} \right)^{\frac{1}{n+1}} \left( \frac{1}{\langle x, u(x)\rangle}\right)^{\frac{1}{n+1}} \nonumber \\
&=& c_n \   |S|_n ^{\frac{2}{n+1}}  \left(\frac{m_{S^\circ} (x)}{m_{S} (x) } \right)^{\frac{1}{n+1}}  = c_n \ \left(|S|_n|S^{\circ}|_n\right)^{\frac{1}{n+1}}  \left(\frac{n_{S^\circ} (x)}{n_{S} (x) } \right)^{\frac{1}{n+1}} ,
\end{eqnarray}
where $m_{S^\circ } $ and  $n_{S^\circ} $  are defined as follows.
Denote by $u_S : \partial S \rightarrow \mathbb{S}^{n-1}$, $x \rightarrow u(x)$ the Gauss map of $S$.  Then, 
$m_{S^\circ} (x) = \frac{1}{n}  \frac{\kappa_S(x) } {\langle x, u(x)\rangle ^n} $ is the density function of the ``cone measure"  $M_{S^\circ}$ of $S^\circ$. 
For a Borel set $A \in \partial S$,   $M_{S^\circ} (A) = | \conv[0, u_{S^\circ}^{-1} (u_S(A))]|_n$ and $n_{S^\circ} (x) = \frac{1}{n |S^\circ|_n}  \frac{\kappa_S(x) } {\langle x, u(x)\rangle^n} $ is the density of the normalized cone measure $\mathbb{P}_{S^{\circ}}$ of $S^{\circ}$ (see e.g. \cite{PaourisWerner2012}). 
This means  that 
$$
d   M_{S^\circ} (x) = m_{S^\circ} (x) d \mu_S(x) \  \   \text {and } \  \    d \mathbb{P}_{S^\circ} (x) = n_{S^\circ} (x) d \mu_S(x).
$$
Thus, the right-hand side of Theorem \ref{main2} can be re-written as an expression involving the normalized cone measures
\begin{eqnarray*} \label{Theorem:smooth}
 \lim _{\delta \rightarrow 0} \frac{\d_S(\delta) -1}{\delta^{\frac{2}{n+1}} }
 = c_n \left(|S|_n|S^{\circ}|_n\right)^{\frac{1}{n+1}} \left[ \max_{x\in\partial S} \left(\frac{n_{S^\circ} (x)}{n_{S} (x) } \right)^{\frac{1}{n+1}}- \min_{x\in\partial S} \left(\frac{n_{S^\circ} (x)}{n_{S} (x) } \right)^{\frac{1}{n+1}}  \right].
 \end{eqnarray*}
\vskip 3mm
\noindent
Theorems  \ref{main1} and \ref {main2} give rise to a new affine invariant which we will also call $G$.  We set
$$
\underline{G}(S)=\liminf _{\delta \rightarrow 0} \frac{\d_S(\delta) -1}{\delta^{\frac{2}{n+1}} } \hskip 3mm  \text{    and   } \hskip 3mm \overline{G}(S)=\limsup _{\delta \rightarrow 0} \frac{\d_S(\delta) -1}{\delta^{\frac{2}{n+1}}}
$$
and we put \(G(S)=\underline{G}(S)\),  if \(\underline{G}(S)=\overline{G}(S)\).
By a theorem of Petty \cite{Petty}, a centrally symmetric convex body \(S\) with \(C^2\)-boundary is an ellipsoid if and only if there is a constant \(\alpha_S\) such that \(\alpha_S=\frac{\kappa(x)^{\frac{1}{n+1}}}{\langle x, u(x)\rangle}\) for every \(x\in\partial S\). Therefore,  \(G(S)=0\)  if and only if $S$ is an ellipsoid. An immediate consequence of this fact and Theorem \ref{main2} is the following corollary.
\begin{cor}  Let $S\subseteq \mathbb{R}^n$  be a centrally symmetric convex body with \(C_+^2\)-boundary. Suppose  there exists a constant \(\delta_0\) such that for all \(0<\delta<\delta_0\) and all \(\delta'>0\) we have that
\[
S_{\delta}=\ii{S}.
\]
Then \(S\) is an ellipsoid.
\end{cor}
This corollary supports the conjecture that equality of the  floating body of $S$ and the polar of the illumination body of the polar $S^\circ$  characterizes ellipsoids. Note that in \cite{SchuettWerner1994, Stancu2009, WernerYe11} similar theorems for the homothety conjecture also make use of Petty's lemma as a crucial step in their proofs.

\vskip 5mm

\section{Background}\label{SectionBackground}

Let \(K\subseteq \mathbb{R}^n\) be a convex body and let \(x\in\partial K\). If \(x\) has a unique outer normal,  we denote it by \(u_K(x)\).  We omit the subscript \(K\), if it is clear what is meant from the context. If the outer normal is well-defined everywhere,  we also denote by \(u:\partial K\rightarrow \mathbb{S}^{n-1}\), $ x \rightarrow u(x)$,  the \textit{Gauss map}. 
\par
Let \(u(x)\) be defined, let \(L\in O(n)\) be a rotation such that \(L(u(x))=-e_n\) and put \(K'=L(K-x)\). Then there is  \(\tau>0\) and a convex function \(f_x:\tau B_2^{n-1}\rightarrow \mathbb{R}_{\geq 0}\) with \(f(0)=0\) such that the boundary of \(K'\) is locally around the origin given by the graph of \(f_x\). We call \(f_x\) a \textit{parametrization} of \(K\) at \(x\). 
We say that  \(y\in \partial K\) \textit{corresponds to} \((z, f_x(z))\) , \(z\in\tau B_2^{n-1}\), or vice versa, \((z,f_x(z))\) \textit{corresponds to} \(y\) if \(L(y-x)=(z,f_x(z))\).
\par
If \(K\) is differentiable at \(0\),  then \(\nabla f_x(0)=0\). If \(f_x\) is \(C^2\) on a neighbourhood of \(0\),  the \textit{principal curvatures} and the \textit{Gauss curvature} of \(K\) at \(x\) are defined as the eigenvalues  and the determinant of the Hessian \(Hf_x\) at \(0\). 
We denote the Gauss curvature at \(x\) by \(\kappa_K(x)\) and we omit the subscript  \(K\) in \(\kappa_K\) in most cases since the convex body involved will  usually be clear from the context. 
These definitions are independent of the choice of \(f_x\) (see also the introduction of \cite{ReisnerSchuettWerner, Schuett:2003} for the definition of curvature for convex bodies). 
\par
Provided they exist, we denote by \(\lambda_1(x)\geq \lambda_2(x)\geq\dots\geq \lambda_{n-1}(x)\geq 0\) the principal curvatures at \(x\in\partial K\) and we put \(\lambdam=\min_{x\in\partial K}\lambda_{n-1}(x)\) and \(\lambdaM=\max_{x\in\partial K}\lambda_1(x)\). We put \(\kappam=\min_{x\in\partial K}\kappa(x)\) and \(\kappaM=\max_{x\in\partial K}\kappa(x)\). 
\vskip 3mm
\begin{defn}\label{Delta} Let \(K\subseteq \mathbb{R}^n\) be a convex body and let \(x\in\partial K\) be a boundary point with unique outer normal \(u(x)\). We denote by \(\Delta_x(\delta)>0\) the unique value such that
\[
|\{y\in K: \langle y, u(x) \rangle\geq \langle x, u(x)\rangle-\Delta_x(\delta)\}|_n=\delta |K|_n\quad.
\]
\end{defn}
\vskip 2mm
\begin{defn} Let \(K\subseteq \mathbb{R}^n\) be a convex body with \(0\in\mathrm{int}K\). Let \(x\in\partial K\) and \(\delta\geq 0\). We denote by \(x_{\delta}\in\partial K_{\delta}\) the unique point  such that \(\{x_{\delta}\}=K_{\delta}\cap [0,x]\), by \(x^{\delta}\in\partial K^{\delta}\) the unique point such that \(\{x\}=K\cap [0,x^{\delta}]\) and by \(\ii[\delta]{x}\in\partial\ii[\delta]{K}\) the unique point such that  \(\{\ii[\delta]{x}\}\in[0,x]\cap \ii[\delta]{K}\).
\end{defn}
\vskip 3mm
\noindent 
We call \(y\in\partial K\) a \textit{touching point} of \(x\) with \(K\) if the line segment \([y,x]\) lies in a support hyperplane of \(K\) at \(y\). The following lemma tells us that \(\conv[K,x]\) depends only on the touching points of \(x\) with \(K\).
We include proofs of the next two lemmas even though they are  probably known.
\par
\begin{lem}\label{ConvexHullBdaryII} Let \(K\subseteq\mathbb{R}^n\) be a convex body and let \(x\in\mathbb{R}^n\backslash K\). Then 
\begin{align}
 \partial \left(\mathrm{conv}[K,x]\right)
\subseteq \partial K \cup \{\lambda x+(1-\lambda)y:\lambda\in[0,1], y\in \partial K \text{ is a touching point} \}\quad.\notag
\end{align}
\end{lem} 
\textit{Proof.} Let \(z\in\partial \left(\mathrm{conv}[K,x]\right)\) and assume that \(z\neq x\) and \(z\not \in \partial K\). It is an elementary fact of convex geometry  that \(\mathrm{conv}[K,x]=\{\lambda x+(1-\lambda)y:y\in K,\lambda\in[0,1]\}\). A  simple consequence  is that \(\partial \left(\mathrm{conv}[K,x]\right)\subseteq \{\lambda x+(1-\lambda)y:y\in\partial K,\lambda\in [0,1]\}\) (see,  e.g.,  \cite{Webster}). Hence there are \(\lambda \in (0,1)\) and \(y\in\partial K\) such that \(z=\lambda y +(1-\lambda) x\). Let \(H\) be a support hyperplane of \(\mathrm{conv}[K,x]\) at \(z\) and let \(H^+\) be the corresponding closed halfspace including \(\mathrm{conv}[K,x]\). Then \(y,x\in H^+\), hence, \(y,x\in H\) and it follows that \(H\) is also a support hyperplane of \(\conv[K,x]\) at \(y\). \(\Box\)
\par
\bigskip
\noindent
\begin{lem}\label{1dHeightestimate}
Let \(\tau>0\) and let \(f:(-\tau,\tau)\rightarrow\mathbb{R}_{\geq 0}\) be convex, of class \(C^2\) and such that \(f(0)=0\), \(f'(0)=0\) and \(f''(0)>0\). Let \(\frac{1}{16}>\eta\geq 0\). Let \(\tau>\theta>0\)
be such that \(|t|\leq \theta\) implies 
$$\frac{1-\eta}{2}f''(0)t^2\leq f(t)\leq \frac{1+\eta}{2}f''(0)t^2.
$$
If \(\Delta > 0\) is such that \(\Delta < \frac{(1-4\sqrt{\eta})f''(0)}{2}\theta^2\), then there exists some \(0 < t_0 <\theta\) such that the line through \(-\Delta e_2\) and \((t_0, f(t_0))\) is  tangent to \(f\) at \((t_0,f(t_0))\). In this case
\(f(t_0)\leq \frac{1+\eta}{1-4\sqrt{\eta}}\Delta\).
\end{lem}
\textit{Proof.} For every \(s\in (0,\tau)\) the line \(l_{s}\) through \(-\Delta e_2\) and \((s,f(s))\) is given by \(l_s(t)=\frac{f(s)+\Delta}{s}t-\Delta\) and \(l_s\) touches the graph of \(f\) if and only if \(f'(s)=\frac{f(s)+\Delta}{s}\). We provide upper and lower estimates for \(f'(t)\). Put \(c=1-2\sqrt{\eta}\). Since \(f\) is convex we have for every \(0<t\leq \theta\) that
\begin{align}
f'(t)\geq &
\frac{f(t)-f(tc)}{t-tc}
\geq \frac{\frac{1-\eta}{2}f''(0)t^2-\frac{1+\eta}{2}f''(0)c^2t^2}{(1-c)t}
=\left(\frac{1+c}{2}-\eta\frac{1+c^2}{2(1-c)}\right)f''(0)t\notag\\
\geq & \left(\frac{1+c}{2}-\frac{\eta}{1-c}\right)f''(0)t
=\left(1-\sqrt{\eta}-\frac{1}{2}\sqrt{\eta}\right)f''(0)t=\left(1-\frac{3}{2}\sqrt{\eta}\right)f''(0)t\quad.\notag
\end{align}
Put \(d=1+2\sqrt{\eta}\). Providing that \(0<(1+2\sqrt{\eta})t\leq \theta\) we get the following upper estimate.
\begin{align}
f'(t)\leq & \frac{f(dt)-f(t)}{dt-t} 
\leq \frac{\frac{1+\eta}{2}f''(0)d^2t^2-\frac{1-\eta}{2}f''(0)t^2}{dt-t}
=\left(\frac{1+d}{2}+\eta\frac{1+d^2}{2(d-1)}\right)f''(0)t\notag\\
=&\left(1+\sqrt{\eta}+\frac{\sqrt{\eta}}{2}+\eta+\eta^{3/2}\right)f''(0)t\leq \left(1+\frac{7}{2}\sqrt{\eta}\right)f''(0)t\notag\quad.
\end{align}
Put \(t=\sqrt{\frac{2\Delta}{f''(0)(1-4\sqrt{\eta})}}\). Then \(t<\theta\) and we get
\begin{align}
f'(t)t-f(t)\geq \left(1-\frac{3}{2}\sqrt{\eta}\right)f''(0)t^2-\frac{1+\eta}{2}f''(0)t^2\geq \left(1-4\sqrt{\eta}\right)\frac{f''(0)}{2}t^2=\Delta\quad.\notag
\end{align}
Now, put \(t=\sqrt{\frac{2\delta}{f''(0)(1+8\sqrt{\eta})}}\). It is easy to verify that \(\frac{t}{1+2\sqrt{\eta}}\leq \theta\). It follows that
\begin{align}
f'(t)t-f(t)\leq \left(1+\frac{7}{2}\sqrt{\eta}\right)f''(0)t^2-\frac{1-\eta}{2}f''(0)t^2\leq (1+8\sqrt{\eta})\frac{f''(0)}{2}t^2=\Delta\quad.\notag
\end{align}
Since \(f'(t)t-f(t)\) is continuous there is \(\sqrt{\frac{2\delta}{f''(0)(1+8\sqrt{\eta})}}\leq t_0\leq \sqrt{\frac{2\Delta}{f''(0)(1-4\sqrt{\eta})}}\) such that \(f'(t_0)=\frac{f(t_0)+\Delta}{t_0}\). 
Note that by monotonicity of \(f\) we obtain
\[
f(t_0)\leq f\left(\sqrt{\frac{2\Delta}{f''(0)(1-4\sqrt{\eta})}}\right)
\leq\frac{1+\eta}{2}f''(0)\sqrt{\frac{2\Delta}{f''(0)(1-4\sqrt{\eta})}}^{\ 2}\frac{1+\eta}{1-4\sqrt{\eta}}\Delta\quad.
\]
\hfill \(\Box\)

\bigskip
\noindent
We obtain the following immediate generalization for higher dimensions:
\begin{cor}\label{MultidimHeightEstimate} Let \(\tau>0\) and let \(f:\mathrm{int}(\tau B_2^{n-1})\rightarrow \mathbb{R}_{\geq 0}\) be convex, of class \(C^2\) such that \(f(0)=0\), \(\nabla f(0)=0\) and the smallest eigenvalue \(\lambda_{n-1}\) of \(Hf(0)\) is positive. Let \(\frac{1}{16}>\eta\geq 0\) and let \(\tau>\theta>0\) be such that \(\|z\|_2\leq \theta\) implies 
$$
\frac{1-\eta}{2}\langle Hf(0)z,z\rangle\leq f(z)\leq \frac{1+\eta}{2}\langle Hf(0)z,z\rangle.
$$ 
If \(\Delta>0\) is such that \(\Delta<\frac{(1-4\sqrt{\eta})\lambda_{n-1}}{2}\theta^2\), then for every \(v\in \mathbb{S}^{n-2}\) there is some \(0<t_0(v)<\theta\) such that the line through \(-\Delta e_n\) and \((t_0(v)v,f(t_0(v)v))\) lies in the tangent hyperplane of \(f\) at \((t_0(v)v,f(t_0(v)v))\). In this case \(f(t_0(v)v)\leq \frac{1+\eta}{1-4\sqrt{\eta}}\Delta\).                                 
\end{cor} 
\bigskip
\noindent
The following lemma can be found in [\cite{Reitzner:2002}, Lemma 6]. 
\begin{lem}\label{OvaloidParam} Let \(K\subseteq \mathbb{R}^n\) a convex body with \(C_+^2\)-boundary. There is  \(\tau>0\) such that for every \(\xi\in \partial K\) there is a parametrization 
\(f_{\xi}: \tau B_2^{n-1}\rightarrow \mathbb{R}_{\geq 0}\) of \(K\) at \(\xi\) such that for every \(1>\eta>0\) there is  \(\tau>\theta>0\),  independent of \(\xi\),  such that \(\|z\|_2\leq \theta\) implies 
$$
\frac{1-\eta}{2}\langle Hf_{\xi}(0)z,z\rangle \leq f_{\xi}(z)\leq \frac{1+\eta}{2}\langle Hf_{\xi}(0)z,z\rangle.
$$
\end{lem}
\vskip 3mm
\noindent
A careful analysis of the proof in \cite{Reitzner:2002} yields the following version of this lemma for \(C^2\)-bodies.
\par
\begin{lem}\label{BodyParam} Let \(K\subseteq \mathbb{R}^n\) be a convex body with \(C^2\) boundary. Then there is  \(\tau>0\) such that for every \(\xi\in \partial K\) there is a parametrization 
\(f_{\xi}: \tau B_2^{n-1}\rightarrow \mathbb{R}_{\geq 0}\) of \(K\) at \(\xi\) such that for every \(\eta>0\) there is a $\theta$ independent of \(\xi\),  \(\tau>\theta>0\),  such that \(\|z\|_2\leq \theta\) implies 
\[
\frac{1}{2}\left(\langle Hf_{\xi}(0)y,y\rangle-\eta\|y\|^2\right)\leq f_{\xi}(z)\leq \frac{1}{2}\left(\langle Hf_{\xi}(0)y,y\rangle+\eta\|y\|^2\right)\quad.
\]
\end{lem}
\vskip 3mm
As eventually we treat symmetric convex bodies with smooth boundary,  we will from now on mostly consider symmetric convex bodies $S$ that are $C^2$,  even though some of the mentioned results hold true for general convex bodies.
\par
\begin{cor}\label{ConvexHullLemma} Suppose that \(S\) has \(C_+^2\)-boundary. For every \(\varepsilon >0\) there is \(\Delta_0>0\) such that for every \(0\leq \Delta\leq \Delta_0\)  and every \(x\in\partial S\) the following holds: Let 
\[
S_1=S\cap \{\xi\in\mathbb{R}^n:\langle \xi-x, u(x)\rangle\leq -(1+\varepsilon)\Delta\}
\]
and 
\[
S_2=S\cap \{\xi\in\mathbb{R}^n:\langle \xi-x, u(x)\rangle\geq -(1+\varepsilon)\Delta\}\quad.
\]
Then
\[
\conv[S, x+\Delta u(x)]=S_1\cup\conv[S_2,x+\Delta u(x)]\quad.
\]
\end{cor}
\textit{Proof.} Let \(\varepsilon>0\) and \(\frac{1}{3}>\eta>0\) with \(\frac{1+\eta}{1-4\sqrt{\eta}}\leq 1+\varepsilon\). 
By Lemma \ref{OvaloidParam}, we can choose \(\tau>0\) and  \(\tau>\theta>0\) such that for every \(x\in \partial S\)  the following holds:
\par
\noindent
There is a parametrization \(f_x:\tau B_2^{n-1}\rightarrow\mathbb{R}_{\geq 0}\) such that  
$$
\frac{1-\eta}{2}\langle Hf_x(0) z,z\rangle\leq f_x(z)\leq \frac{1+\eta}{2}\langle Hf_x(0)z,z\rangle, 
$$
if \(\|z\|_2\leq \theta\). Put \(\Delta_0=\frac{(1-4\sqrt{\eta})\lambdam\theta}{2}\). 
Let \(x_0\in\partial S\) and \(0<\Delta<\Delta_0\). By Lemma \ref{ConvexHullBdaryII} it is sufficient to show that the touching points of \(x_0+\Delta u(x_0)\) with \(S\) lie in \(S_2\). If \(y\in \partial S\) is a touching point,  then there is some \(w\in\mathbb{S}^{n-1}\) orthonormal to \(u(x_0)\) such that 
\[
y\in \partial S\cap \{x_0+\mu_1u(x_0)+\mu_2w:\mu_1,\mu_2\in\mathbb{R},\mu_2\geq 0\}\quad.
\]
It is obvious that \(y\) is the unique touching point which lies on the halfplane \(\{x_0+\mu_1u(x_0)+\mu_2w:\mu_1,\mu_2\in\mathbb{R},\mu_2\geq 0\}\). Let \(f_{x_0}:\tau B_2^{n-1}\rightarrow \mathbb{R}_{\geq 0}\) be a parametrization of the boundary at \(x_0\) and let \(v\in\mathbb{S}^{n-2}\) be the vector corresponding to \(w\). By Corollary \ref{MultidimHeightEstimate} there is a \(t_0(v)\) such that the line through \(-\Delta e_n\) and \((t_0(v)v, f_{x_0}(t_0(v)))\) is tangential to the graph of \(f_{x_0}\). Hence, \((t_0(v)v,f_{x_0}(t_0(v)))\) corresponds to the touching point \(y\) and \(f_{x_0}(t_0(v)v)\leq \frac{1+\eta}{1-4\sqrt{\eta}}\Delta\leq (1+\varepsilon)\Delta\). It follows that
\[
y\in \{\xi\in\mathbb{R}^n:\langle\xi, u(x_0)\rangle\geq \langle x_0, u(x_0)\rangle-(1+\varepsilon)\Delta\}\quad.
\]  
\hfill \(\Box\)
\par
\bigskip
\noindent
We use Corollary \ref{ConvexHullLemma} to obtain an upper volume estimate for the convex hull of \(S\) with a point. Let \(\tau=\tau(S)>0\) and \(f_x:\tau B_2^{n-1}\rightarrow \mathbb{R}_{\geq 0}\) be chosen according to Lemma \ref{OvaloidParam}. We may assume without loss of generality that \(\tau\) is chosen so small such that for every \(z\in\tau B_2^{n-1}\) it holds that \(\frac{1}{4}\langle Hf_x(0)z,z\rangle\leq f_x(z)\leq \frac{3}{4}\langle Hf_x(0)z,z\rangle\).
For every \(v\in\mathbb{S}^{n-2}\),
\begin{align}
f_{x}(\tau v)\geq \frac{1}{4}\langle Hf_x(0)\tau v,\tau v\rangle\geq\frac{\lambdam \tau^2}{4}=:T_0(S)\quad.\notag
\end{align}
Hence, the part of the boundary of \(S\) lying in the halfspace \(\{\xi\in\mathbb{R}^n: \langle \xi,u(x)\rangle\geq \langle x, u(x)\rangle -T_0(S)\}\)
is completely parametrized by \(f_x\), i.e.,  for every 
\[
y\in\partial S\cap \{\xi\in\mathbb{R}^n: \langle \xi,u(x)\rangle\geq \langle x, u(x)\rangle -T_0(S)\}
\]
there is a \(z\in \tau B_2^{n-1}\) such that \((z, f_x(z))\) corresponds to \(y\). 
\vskip 2mm
\noindent
The following two lemmas are uniform versions of well-known estimates for cap and hat volumes which can be found in [\cite{Leichtweiss}, p. 459] and [\cite{Werner94}, Lemma 3]. 
\begin{lem}\label{CapEstimates} There is a non-negative function \(\phi\) with \(\lim_{\Delta\rightarrow 0}\phi(\Delta)=0\) such that for every \(x\in\partial S\)
\[
\frac{(2\Delta)^{\frac{n+1}{2}}}{n+1}\frac{|B_2^{n-1}|_{n-1}}{(\kappa(x)+\phi(\Delta))^{\frac{1}{2}}}
\leq|\{\xi\in S: \langle \xi-x, u(x)\rangle\geq -\Delta\}|_n\quad.
\]
\end{lem}
\textit{Proof.} We show that for every \(\varepsilon>0\) there is a \(\Delta_0>0\) such that for every \(0\leq \Delta\leq\Delta_0\) we have
\[
\frac{(2\Delta)^{\frac{n+1}{2}}}{n+1}\frac{|B_2^{n-1}|_{n-1}}{(\kappa(x)+\varepsilon)^{\frac{1}{2}}}
\leq|\{\xi\in S: \langle \xi-x, u(x)\rangle\geq -\Delta\}|_n
\]
which establishes the proof. Put \(\varepsilon'=\frac{\varepsilon}{\lambdaM 2^{n-1}}\) and note that \(\lambdaM>0\). Assume without loss of generality that \(\varepsilon>0\) is chosen sufficiently small so that \(\frac{\varepsilon'}{\lambdaM}\leq 1\). Apply Lemma \ref{BodyParam} to \(\eta=\frac{\varepsilon'}{2}\) and let \(\tau>\theta>0\) be chosen accordingly to this \(\eta\). Let \(Hf_x(0)=\sum_{i=1}^{n-1}\lambda_i(x)v_i\otimes v_i\). Put \(\tilde{\lambda_i}(x)=\max[\lambda_i(x),\frac{\varepsilon'}{2}]\) and  \(\tilde{H}_x=\sum_{i=1}^{n-1}\tilde{\lambda_i}(x)v_i\otimes v_i\). For \(w\in\mathbb{S}^{n-2}\) we get
\begin{eqnarray*}
&&\hskip -5mm \{s\in [-\tau,\tau]: f_x(sw)\leq t\}\supseteq \{s\in [-\theta,\theta]:f_x(sw)\leq t\}
\supseteq \left\{s\in [-\theta,\theta]:\frac{s^2}{2}\left(\langle Hf_x(0)w,w\rangle+\frac{\varepsilon'}{2}\right)\leq t\right\}\\
&&\hskip 3mm \supseteq \left\{s\in [-\theta,\theta]:\frac{s^2}{2}\left(\langle \tilde{H}_x w,w\rangle+\frac{\varepsilon'}{2}\right)\leq t \right\}
= \left\lbrace s\in [-\theta,\theta]: |s|\leq \sqrt{\frac{2t}{\langle \tilde{H}_x w,w \rangle+\frac{\varepsilon}{2}}}\right\rbrace .
\end{eqnarray*}
The second inclusion follows as 
\[
f_x(sw)\leq \frac{1}{2}\left(\langle Hf_x(0)sw,sw\rangle+\frac{\varepsilon'}{2}\|sw\|^2\right)=\frac{s^2}{2}\left(\langle Hf_x(0)w,w\rangle+\frac{\varepsilon'}{2}\right)\quad.
\]
We may assume that \(t\leq \frac{\varepsilon'}{4}\theta^2\). Note that \(s\leq \theta\) since \({\langle \tilde{H}_x w,w \rangle+\frac{\varepsilon'}{2}}\leq \frac{\varepsilon'}{2}\|w\|_2^2+\frac{\varepsilon'}{2}=\varepsilon'\). This yields
\[
\{z\in \tau B_2^{n-1}: f_x(z)\leq t\}\supseteq \sqrt{2t} \ll y\in\mathbb{R}^{n-1}: \left\langle \left(\tilde{H}_x+\frac{\varepsilon'}{2}I\right)y,y \right\rangle\leq 1\rr
\] 
and the set on the right-hand side is an ellipsoid with principle axes 
\(
\left(\tilde{\lambda}_i(x)+\frac{\varepsilon'}{2}\right)^{-1/2}v_i
\). Note that 
\(
\left(\tilde{\lambda}_i(x)+\frac{\varepsilon'}{2}\right)^{-1/2}\geq (\lambda_i(x)+\varepsilon')^{-1/2}
\),
which, for \(\Delta\leq \frac{\varepsilon'}{4}\theta^2\), implies  that
\begin{eqnarray*}
\int_0^{\Delta} |\{z\in\tau B_2^{n-1}:f_x(z)\leq t\}|_{n-1}\mathrm{d}t
&\geq&  \int_0^{\Delta}(2t)^{\frac{n-1}{2}}\left|\left\lbrace y\in\mathbb{R}^{n-1}: \left\langle \left(\tilde{H}_x+\frac{\varepsilon'}{2}I\right)y,y \right\rangle\leq 1\right\rbrace\right|_{n-1}\mathrm{d}t\\
&\geq& \frac{(2\Delta)^{\frac{n+1}{2}}}{n+1}\prod_{i=1}^{n-1}\frac{1}{(\lambda_i(x)+\varepsilon')^{1/2}}|B_2^{n-1}|_{n-1}\quad.
\end{eqnarray*}
We can conclude, as
\begin{eqnarray*}
\prod_{i=1}^{n-1}(\lambda_i(x)+\varepsilon') &\leq& \prod_{i=1}^{n-1}\lambda_i(x)+\prod_{i=1}^{n-1}(\lambdaM+\varepsilon')-\lambdaM^{n-1}\\
&=&\kappa(x)+\lambdaM^{n-1}\left(1+\frac{\varepsilon'}{\lambdaM}\right)^{n-1}-\lambdaM^{n-1}
=\kappa(x)+\lambdaM^{n-1}\sum_{i=1}^{n-1}\binom{n-1}{i}\left(\frac{\varepsilon'}{\lambdaM}\right)^i \\
&\leq&\kappa(x)+\lambdaM^{n-1}\sum_{i=1}^{n-1}\binom{n-1}{i}\frac{\varepsilon'}{\lambdaM}
\leq\kappa(x)+\varepsilon'\lambdaM^{n-2}\sum_{i=0}^{n-1}\binom{n-1}{i}\\
&=&\kappa(x)+\varepsilon'\lambdaM^{n-2}2^{n-1}=\kappa(x)+\varepsilon\quad.
\end{eqnarray*}
\hfill \(\Box\)

\begin{lem}\label{CapHatEstimates} Suppose that \(S\) has \(C_+^2\)-boundary. Then there are non-negative functions \(\phi\) and \( \psi\) with \(\lim_{\Delta\rightarrow 0}\phi(\Delta)=0\) and \(\lim_{\Delta\rightarrow 0}\psi(\Delta)=0\) such that for every \(x\in \partial S\) the following holds.
\begin{enumerate}
\item 
\begin{eqnarray*}
\hskip -4mm \frac{(2\Delta)^{\frac{n+1}{2}}}{n+1}\frac{|B_2^{n-1}|_{n-1}}{\kappa(x)^{\frac{1}{2}}}(1-\phi(\Delta))
\leq|\{\xi\in S: \langle \xi-x, u(x)\rangle\geq -\Delta\}|_n
\leq \frac{(2\Delta)^{\frac{n+1}{2}}}{n+1}\frac{|B_2^{n-1}|_{n-1}}{\kappa(x)^{\frac{1}{2}}}(1+\phi(\Delta)).
\end{eqnarray*}
\item 
For every \(y\in \mathbb{R}^n\) such that \(\langle y, u(x)\rangle\geq \langle x, u(x)\rangle+\Delta\) we have
\begin{equation*}
\frac{(2\Delta)^{\frac{n+1}{2}}}{n(n+1)}\frac{|B_2^{n-1}|_{n-1}}{\kappa(x)^{\frac{1}{2}}}(1-\psi(\Delta))\leq|\conv[S,y]|_n-|S|_n\quad.
\end{equation*}
\item 
\begin{align}
|\conv[S,x+\Delta u(x)]|_n-|S|_n
\leq  \frac{(2\Delta)^{\frac{n+1}{2}}}{n(n+1)}\frac{|B_2^{n-1}|_{n-1}}{\kappa(x)^{\frac{1}{2}}}(1+\psi(\Delta))\notag .
\end{align}
\end{enumerate}
\end{lem}
\par
\noindent
\textit{Proof.} Before we show the volume estimates,  we establish some general facts.  
\par
\noindent
Let \(0<\varepsilon<1\) and \(\Delta_0>0\) be chosen according to Corollary \ref{ConvexHullLemma}. Let \(0<\theta<\tau\) be such that for every \(z\in\theta B_2^{n-1}\) and every \(x\in\partial S\) we have 
$$
\frac{1-\varepsilon}{2}\langle Hf_x(0)z,z\rangle\leq f_x(z)\leq \frac{1+\varepsilon}{2}\langle Hf_x(0)z,z\rangle .
$$
Without loss of generality we can assume  that  \((1+\varepsilon)\Delta_0\leq T_0(S)\) and \(\Delta_0\leq \frac{(1-\varepsilon)\lambdam \theta^2}{2}\). Let \(0\leq t\leq (1+\varepsilon)\Delta_0\). We show that for every \(x\in\partial S\), 
\begin{align}
\sqrt{\frac{2t}{1+\varepsilon}}\mathcal{E}_x\subseteq \{z\in\tau B_2^{n-1}:f_x(z)\leq t\}\subseteq \sqrt{\frac{2t}{1-\varepsilon}}\mathcal{E}_x , \label{IndicatrixEstimate}
\end{align}
where \(\mathcal{E}_x=\{\zeta\in\mathbb{R}^{n-1}:\langle Hf_x(0)\zeta,\zeta\rangle\leq 1\}\) is the indicatrix of Dupin at \(x\). 
For \(z\in \tau B_2^{n-1}\backslash \{0\}\) it follows for \(z\in \theta B_2^{n-1}\) that
\begin{eqnarray*}
\frac{(1-\varepsilon)\lambdam\|z\|_2^2}{2} &\leq&\frac{1-\varepsilon}{2}\langle Hf_x(0)z,z \rangle\leq f_x(z)\leq \frac{1+\varepsilon}{2} \langle Hf_x(0)z, z\rangle \\
&\leq& \frac{(1+\varepsilon)\lambdaM \|z\|_2^2}{2}\quad.
\end{eqnarray*}
Since \(\sigma\mapsto f_x(\sigma v)\) is strictly monotonously increasing on \([0,\tau]\) for \(v\in \mathbb{S}^{n-2}\),  we conclude for  \(\|z\|_2>\theta\) and \(z'=\frac{z}{\|z\|_2}\) that
\[
f_x(z)>f_x(\theta z')\geq \frac{(1-\varepsilon)\lambdam \theta^2}{2}\geq \Delta_0\geq t.
\]
Therefore,
\begin{eqnarray*}
\{z\in\tau B_2^{n-1}: f_x(z)\leq t\}= \{z\in \theta B_2^{n-1}: f_x(z)\leq t\}
\subseteq &\{z\in \theta B_2^{n-1}: \frac{1-\varepsilon}{2}\langle Hf_x(0)z,z\rangle\leq t\} .
\end{eqnarray*}
Since for \(\|z\|_2> \theta\)  we have \(\frac{1-\varepsilon}{2}\langle Hf_x(0)z,z\rangle> \frac{(1-\varepsilon)\lambdam\theta^2}{2}\geq t\),  it follows that
\begin{eqnarray*}
 \ll z\in \theta B_2^{n-1}: \frac{1-\varepsilon}{2}\langle Hf_x(0)z,z\rangle\leq t\rr=\ll z\in \mathbb{R}^{n-1}: \frac{1-\varepsilon}{2}\langle Hf_x(0)z,z\rangle\leq t\rr
=\sqrt{\frac{2t}{1-\varepsilon}}\mathcal{E}_x\quad.
\end{eqnarray*}
We conclude in a similar way that
\(
\{z\in\tau B_2^{n-1}: f_x(z)\leq t\}\supseteq \sqrt{\frac{2t}{1+\varepsilon}}\mathcal{E}_x
\).
\vskip 2mm
\noindent
\textit{Proof of 1.} Let \(x\in\partial S\).  Since \(\Delta_0\leq T_0(S)\) the part of the boundary of \(S\) lying in the halfspace \(\{\xi\in\mathbb{R}^n:\langle \xi-x, u(x)\rangle\geq -\Delta_0\}\) is completely parametrized by \(f_x\).  Hence, for every \(0\leq\Delta\leq \Delta_0\),  the volume of \(\{\xi\in S: \langle \xi-x, u(x) \rangle\geq -\Delta\}\) equals the volume of 
\[
\{(z,t)\in\mathbb{R}^n: f_x(z)\leq t\leq \Delta\}\quad.
\]
Cavalieri's principle and the right-hand side of (\ref{IndicatrixEstimate}) yield
\begin{eqnarray*}
|\{(z,t)\in\mathbb{R}^n: f_x(z)\leq t\leq \Delta\}|_n
&=&\int_0^{\Delta}|\{z\in\tau B_2^{n-1}: f_x(z)\leq t\}|_{n-1}\mathrm{d}t
\leq \int_0^{\Delta}\left|\sqrt{\frac{2t}{1-\varepsilon}}\mathcal{E}_x\right|_{n-1}\mathrm{d}t \\
&=&\int_0^{\Delta}(2t)^{\frac{n-1}{2}}\mathrm{d}t\frac{|\mathcal{E}_x|_{n-1}}{(1-\varepsilon)^{\frac{n-1}{2}}}
=\frac{(2\Delta)^{\frac{n+1}{2}}|\mathcal{E}_x|_{n-1}}{n+1}(1-\varepsilon)^{-\frac{n-1}{2}}\quad.
\end{eqnarray*}
Similarly, using the left-hand side of (\ref{IndicatrixEstimate}), 
one has
\[
|\{(z,t)\in\mathbb{R}^n: f_x(z)\leq t\leq \Delta\}|_n\geq \frac{(2\Delta)^{\frac{n+1}{2}}|\mathcal{E}_x|_{n-1}}{n+1}(1+\varepsilon)^{-\frac{n-1}{2}}\quad.
\]
Since \(\mathrm{det}(Hf_x(0))=\kappa(x)\),  it follows that \(|\mathcal{E}_x|_{n-1}=\kappa(x)^{-1/2}|B_2^{n-1}|_{n-1}\).  This shows the first part of the lemma.
\vskip 2mm
\noindent
\textit{Proof of 2.} Let \(0\leq \Delta\leq \Delta_0\), \(y\in \{\xi\in\mathbb{R}^n\langle\xi-x, u(x)\rangle\geq \Delta\}\) and put \(B_x(\Delta)=\{\xi\in S: \langle \xi-x,u(x) \rangle=-\Delta \}\). The set \(\conv[S,y]\backslash S\) includes the set 
\[
\conv[B_x(\Delta),y]\backslash\{\xi\in S:\langle\ \xi-x, u(x)\rangle\geq -\Delta\}\quad.
\] 
The height of the cone \(\conv[B_x(\Delta),y]\)  is at least \(2\Delta\). Using Part 1. of the lemma, a lower estimate for the volume of \(\conv[S,y]\backslash S\) is
\begin{eqnarray*}
&& \hskip -10mm |\conv[B_x(\Delta),y]|_n-|\{\xi\in S:\langle\ \xi-x, u(x)\rangle\geq -\Delta\}|_n \\
&& \hskip 10mm \geq  \frac{2\Delta|B_x(\Delta)|_{n-1}}{n}-\frac{(2\Delta)^{\frac{n+1}{2}}}{n+1}\frac{|B_2^{n-1}|_{n-1}}{\kappa(x)^{\frac{1}{2}}}(1+\phi(\Delta))\quad.
\end{eqnarray*}
Since \(\Delta_0\leq T_0(S)\),  we get  similar to  the proof of Part 1. that the volume of \(|B_x(\Delta)|_{n-1}\) equals
\(
|\{z\in\tau B_2^{n-1}: f_x(z)\leq \Delta\}|_{n-1}
\). Using the left-hand side of (\ref{IndicatrixEstimate}), we obtain 
\[
|B_x(\Delta)|_{n-1}\geq \frac{(2\Delta)^{\frac{n-1}{2}}|\mathcal{E}_x|_{n-1}}{(1+\varepsilon)^{\frac{n-1}{2}}}=\frac{(2\Delta)^{\frac{n-1}{2}}|B_2^{n-1}|_{n-1}}{\kappa(x)^{1/2}}(1+\varepsilon)^{-\frac{n-1}{2}}, 
\]
which yields
\begin{eqnarray*}
&& \hskip -10mm
\frac{2\Delta|B_x(\Delta)|_{n-1}}{n}-\frac{(2\Delta)^{\frac{n+1}{2}}}{n+1}\frac{|B_2^{n-1}|_{n-1}}{\kappa(x)^{\frac{1}{2}}}(1+\phi(\Delta))\\
&&\hskip 10mm \geq  \frac{(2\Delta)^{\frac{2+1}{2}}|B_2^{n-1}|_{n-1}}{\kappa(x)^{1/2}}
\left(\frac{(1+\varepsilon)^{-\frac{n-1}{2}}}{n}-\frac{1+\phi(\Delta)}{n+1}\right)\quad.
\end{eqnarray*}
The expression 
\(
\frac{(1+\varepsilon)^{-\frac{n-1}{2}}}{n}-\frac{1+\phi(\Delta)}{n+1}
\)
is arbitrarily close to \(\frac{1}{n(n+1)}\),  if \(\Delta_0>0\) is small enough to guarantee that \(\phi(\Delta)\) is sufficiently small and it is possible to choose \(\varepsilon>0\) sufficiently small.
\vskip 2mm 
\noindent
\textit{Proof of 3.} Let \(x\in\partial S\).  Since \(\Delta_0\) is chosen according to Corollary \ref{ConvexHullLemma},  we have that  for every \(0< \Delta\leq\Delta_0\) the volume of \(\conv[S,x+\Delta u(x)]\backslash S\) is the same as the volume of \(\conv[S_2,x+\Delta u(x)]\backslash S_2\),   where 
\[
S_2=\{\xi\in S:\langle \xi-x, u(x)\rangle\geq -(1+\varepsilon)\Delta\}\quad.
\]
Since \((1+\varepsilon)\Delta_0\leq T_0(S)\),  the volume of \(\conv[S_2,x+\Delta u(x)]\) is given by the volume of \(\conv[S_2', -\Delta e_n]\), 
where 
\[
S_2'=\{(z,t)\in\mathbb{R}^{n-1}\times \mathbb{R}: f_x(z)\leq t\leq (1+\varepsilon)\Delta\}\quad.
\]
Let \(F_x(\Delta)\subseteq \mathbb{R}^n\) be the cone with base \(B_x(\Delta)=\sqrt{\frac{\Delta}{2(1-\varepsilon)}}(2+\varepsilon)\mathcal{E}_x\times\{(1+\varepsilon)\Delta\}\) and apex \(-\Delta e_n\), i.e.,  \(F_x(\Delta)=\conv[B_x(\Delta),-\Delta e_n]\). It follows from the right-hand side of (\ref{IndicatrixEstimate}) that for every \(0\leq t\leq (1+\varepsilon)\Delta\), 
\[
\{z\in\mathbb{R}^{n-1}:(z,t) \in S_2'\}\subseteq \{z\in\mathbb{R}^{n-1}:(z,t) \in F_x(\Delta)\}.
\] 
Thus, \(S_2'\subseteq F_x(\Delta)\) and \(\conv[S_2', -\Delta e_n]\subseteq F_x(\Delta)\). Since the height of \(F_x(\Delta)\) is \((2+\varepsilon)\Delta\),  the volume of \(F_x(\Delta)\) is given by 
\[
\frac{1}{n}(2+\varepsilon)\left[\sqrt{\frac{2\Delta}{(1-\varepsilon)}}(2+\varepsilon)\right]^{n-1}|\mathcal{E}_x|_{n-1}=\frac{ (2\Delta)^{\frac{n+1}{2}}|B_2^{n-1}|}{\kappa(x)^{1/2}}
\frac{(1+\frac{\varepsilon}{2})^{n}}{n(1-\varepsilon)^{\frac{n-1}{2}}}\quad.
\]
By Part 1.,  the volume of \(S_2\) can be bounded from below by
\[
\frac{(2\Delta)^{\frac{n+1}{2}}|B_2^{n-1}|_{n-1}}{\kappa(x)^{\frac{1}{2}}}\frac{(1-\phi((1+\varepsilon)\Delta))(1+\varepsilon)^{\frac{n+1}{2}}}{n+1}\quad.
\]
Similar to the proof of Part 2.,  one can derive the claim of the third part. 

\hfill \(\Box\)

\section{Upper Bound} \label{UpperBound}
Let \(\Delta_x(\delta)\) be as in Definition \ref{Delta}.
\begin{lem}\label{FloatingBodyElemLowerUpperBound} Let \(\Gamma_{\min}(\delta)=\min_{x\in\partial S}\frac{\Delta_x(\delta)}{\langle x, u(x)\rangle}\) and \(\Gamma_{\max}(\delta)=\max_{x\in\partial S}\frac{\Delta_x(\delta)}{\langle x,u(x)\rangle}\). Then
\begin{align}
(1-\Gamma_{\max}(\delta))S\subseteq S_{\delta}\subseteq (1-\Gamma_{\min}(\delta))S\quad.\label{ElemInclusionsFloatingBody}
\end{align}
\end{lem}
\textit{Proof.} Let \(x\in \partial S\) and \(\alpha\geq 0\) such that \(0\in \{y\in\mathbb{R}^n: \langle y-x, u\rangle \leq-\alpha\}\).
We show that
\begin{align}
\left(1-\frac{\alpha}{\langle x, u\rangle}\right)S\subseteq  \{y\in\mathbb{R}^n:\langle y-x, u\rangle\leq -\alpha\}\label{InclusionHalfspace}\quad.
\end{align}
Let \(\xi\in S\). Then
\(\langle\xi-x, u\rangle\leq 0\) and \(1-\frac{\alpha}{\langle x, u\rangle}\geq 0\),  since \(0\in x+\{y\in\mathbb{R}^n: \langle y, u\rangle \leq-\alpha\}\). It follows that
\begin{align}
\left\langle\left(1-\frac{\alpha}{\langle x, u\rangle}\right)\xi-x, u\right\rangle 
= \left(1-\frac{\alpha}{\langle x, u\rangle}\right)\langle\xi -x, u\rangle-\alpha\leq -\alpha\notag\quad.
\end{align}
\par
\noindent
Since
\[
S_{\delta}=\bigcap_{x\in\partial S} \{y\in\mathbb{R}^n:\langle y-x, u(x)\rangle\leq -\Delta_x(\delta)\}, 
\]
the left-hand side inclusion of (\ref{ElemInclusionsFloatingBody}) follows immediately from (\ref{InclusionHalfspace}). For the right-hand side inclusion of (\ref{ElemInclusionsFloatingBody}) note that 
\(
x_{\delta}\in\{\xi\in S:\langle\xi-x,u(x)\rangle\leq -\Delta_x(\delta)\}
\).
It follows that
\[
\frac{\|x_{\delta}\|_2}{\|x\|_2}\leq \frac{\langle x_{\delta}, u(x)\rangle}{\langle x, u(x)\rangle}
= \frac{\langle x-\Delta(x,\delta)u(x), u(x)\rangle}{\langle x, (u(x))\rangle}=1-\frac{\Delta(x,\delta)}{\langle x, u(x)\rangle}\leq 1-\Gamma_{\min}(\delta)\quad,
\]
i.e. \(S_{\delta}\subseteq(1-\Gamma_{\min}(\delta))S\)
\hfill \(\Box\)
\hskip 3mm
\noindent
\begin{prop}\label{LowerBoundFloatingBodyC2Case} There is a non-negative function \(\Phi\) with \(\lim_{\delta\rightarrow 0}\Phi(\delta)=0\) such that 
\[
\left(1-(1+\Phi(\delta))\GM\delta^{\frac{2}{n+1}}\right)S\subseteq S_{\delta}\quad.
\]
\end{prop}
\textit{Proof.}
We show that for every \(\varepsilon>0\) there is a \(\delta_0>0\) such that for every \(0\leq\delta\leq\delta_0\), 
\[
\left(1-(1+\varepsilon)\GM\delta^{\frac{2}{n+1}}\right)S\subseteq S_{\delta}\quad.
\]  
Let \(x_{\max}\in\partial S\) be such that  \(G(x_{\max})=\GM\). Then \(\kappa(x_{\max})>0\). Put \(\mu=\min_{x\in\partial S}\langle x, u(x)\rangle > 0\) and \(\varepsilon'=\frac{\mu^{n+1}\kappa(x_{\max})\varepsilon}{\langle x_{\max}, u(x_{\max})\rangle^{n+1}}\).
By Lemma \ref{CapEstimates} there is a \(\Delta_0\) such that for every \(0\leq\Delta\leq \Delta_0\) and every \(x\in\partial S\), 
\[
\frac{(2\Delta)^{\frac{n+1}{2}}}{n+1}\frac{|B_2^{n-1}|_{n-1}}{(\kappa(x)+\varepsilon')^{\frac{1}{2}}}
\leq|\{\xi\in S: \langle \xi-x, u(x)\rangle\geq -\Delta\}|_n\quad.
\]
Choose \(\delta_0\) such that for every \(0\leq\delta\leq\delta_0\) and every \(x\in\partial S\) we have  that \(\Delta_x(\delta)\leq \Delta_0\).   Hence, 
\[
\delta|S|_n\geq \frac{(2\Delta_x(\delta))^{\frac{n+1}{2}}}{n+1}\frac{|B_2^{n-1}|_{n-1}}{(\kappa(x)+\varepsilon')^{\frac{1}{2}}}\quad.
\]
This yields
\begin{align}
\frac{\Delta_x(\delta)}{\langle x, u(x)\rangle}\leq &c(S,n)\frac{(\kappa(x)+\varepsilon')^{\frac{1}{n+1}}}{\langle x,u(x)\rangle}\delta^{\frac{2}{n+1}}
= c(S,n)\left(\frac{\kappa(x)}{\langle x, u(x)\rangle^{n+1}}+\frac{\varepsilon'}{\langle x, u(x)\rangle^{n+1}}\right)^{\frac{1}{n+1}}\delta^{\frac{2}{n+1}}\notag\\
\leq& c(S,n)\left(\frac{\kappa(x)}{\langle x, u(x)\rangle^{n+1}}+\frac{\varepsilon\kappa(x_{\max})}{\langle x_{\max},u(x_{\max})\rangle^{n+1}}\right)^{\frac{1}{n+1}}\delta^{\frac{2}{n+1}}\notag\\
\leq& c(S,n)\left(\frac{\kappa(x_{\max})}{\langle x_{\max}, u(x_{\max})\rangle^{n+1}}+\frac{\varepsilon\kappa(x_{\max})}{\langle x_{\max},u(x_{\max})\rangle^{n+1}}\right)^{\frac{1}{n+1}}\delta^{\frac{2}{n+1}}\notag\\
= & c(S,n)\GM(1+\varepsilon)^{\frac{1}{n+1}}\delta^{\frac{2}{n+1}}\leq c(S,n)\GM(1+\varepsilon)\delta^{\frac{2}{n+1}}\quad.\notag
\end{align}
We conclude with  Lemma \ref{FloatingBodyElemLowerUpperBound}.
\hfill \(\Box\)
\vskip 3mm
\noindent
\begin{prop}\label{FloatingBodyOvaloid} Suppose that \(S\) has \(C_+^2\)-boundary. Then there is a non-negative function \(\Phi\) with \(\lim_{\delta\rightarrow 0}\Phi(\delta)=0\) such that
\begin{align}
\left(1-(1+\Phi(\delta))\GM \delta^{\frac{2}{n+1}}\right)S\subseteq S_{\delta}\subseteq \left(1-(1-\Phi(\delta))\Gm\delta^{\frac{2}{n+1}}\right)S . \label{FloatingOvaloidEstimate}
\end{align}
\end{prop} 
\textit{Proof.} Let \(1>\varepsilon>0\) and let \(\phi\) be the function of Lemma \ref{CapHatEstimates}.
Let \(\Delta_0=\Delta_0(\varepsilon)>0\) be sufficiently small, such that for every \(0\leq\Delta\leq \Delta_0\) we have \(\phi(\Delta)\leq \varepsilon\).  Let \(\delta_0>0\) be such that
\[
\frac{(2\Delta_0)^{\frac{n+1}{2}}}{n+1}\frac{|B_2^{n-1}|_{n-1}}{\kappaM^{1/2}}(1-\varepsilon)=\delta_0 |S|_n\quad.
\]
Let \(x\in\partial S\). If \(\delta\leq \delta_0\), then \(\Delta_x(\delta)\leq \Delta_0\). Indeed, if \(\Delta_x(\delta)>\Delta_0\) then 
\begin{eqnarray*}
|\{\xi\in S: \langle \xi-x, u(x)\rangle\geq -\Delta_x(\delta)\}|_n
&>& |\{\xi\in S: \langle \xi-x, u(x)\rangle\geq -\Delta_0\}|_n \\
&\geq& \frac{(2\Delta_0)^{\frac{n+1}{2}}}{n+1}\frac{|B_2^{n-1}|_{n-1}}{\kappaM^{1/2}}(1-\varepsilon)=\delta_0 |S|_n\quad.
\end{eqnarray*}
Since for \(\delta\leq \delta_0\) we have \(\Delta_x(\delta)\leq\Delta_0\),  we obtain the estimate
\begin{align}
\frac{\Gm}{(1+\varepsilon)^{\frac{2}{n+1}}}\delta^{\frac{2}{n+1}}
\leq \frac{c(S,n)}{(1+\varepsilon)^{\frac{2}{n+1}}}\cdot\frac{\kappa(x)^{\frac{1}{n+1}}}{\langle x, u(x)\rangle}\delta^{\frac{2}{n+1}}
\leq\frac{\Delta_x(\delta)}{\langle x, u(x)\rangle}
\leq \frac{c(S,n)}{(1-\varepsilon)^{\frac{2}{n+1}}}\cdot\frac{\kappa(x)^{\frac{1}{n+1}}}{\langle x, u(x)\rangle}\delta^{\frac{2}{n+1}}
\leq \frac{\GM}{(1-\varepsilon)^{\frac{2}{n+1}}}\delta^{\frac{2}{n+1}}\notag
\end{align}
Now apply Lemma \ref{FloatingBodyElemLowerUpperBound}. 
\hfill \(\Box\)

\bigskip
\noindent
\begin{prop}\label{UpperLowerEstimateIlluBody} Suppose that \(S\) has \(C_+^2\)-boundary. Then there is a non-negative function \(\tilde{\Psi}\) with \(\lim_{\delta\rightarrow 0} \tilde{\Psi}(\delta)=0\) such that
\begin{align}
\left(1+(1-\tilde{\Psi}(\delta))n^{\frac{2}{n+1}}\Gm\delta^{\frac{2}{n+1}}\right)S\subseteq  S^{\delta}\subseteq  \left(1+(1+\tilde{\Psi}(\delta))n^{\frac{2}{n+1}}\GM\delta^{\frac{2}{n+1}}\right)S\label{IlluOvaloidEstimate}
\end{align}
\end{prop}
\textit{Proof.} For \(x\in\partial S\), let \(\Delta^x(\delta)\geq 0\) be defined as the value such that
\[
|\conv[S,x+\Delta^x(\delta)u(x)]|_n-|S|_n=\delta|S|_{n}\quad.
\]
Let \(\varepsilon>0\) be given and let \(\Delta_0>0\) be such that for every \(0\leq \Delta\leq \Delta_0\) the function \(\psi\) of Lemma \ref{CapHatEstimates} is smaller than or equal to \(\varepsilon\). 
We show that there is  \(\delta_0\) such that for every \(0\leq \delta\leq \delta_0\) and every  \(y\in\mathbb{R}^n\) with \(|\conv[S,y]|_n-|S|_n=\delta|S|_n\),  it follows that \(\langle y-x, u(x)\rangle\leq \Delta_0\) for every \(x\in\partial S\). In particular, \(\Delta^x(\delta)\leq \Delta_0\). Let \(\delta_0>0\) be such that
\[
\frac{(2\Delta_0)^{\frac{n+1}{2}}}{n(n+1)}\frac{|B_2^{n-1}|_{n-1}}{\kappaM^{1/2}}(1-\varepsilon)=\delta_0 |S|_n\quad.
\]
Then arguments similar to the ones in the proof of Proposition \ref{FloatingBodyOvaloid} ensure that \(\delta_0\) has the desired properties. 
\par
\noindent
Let \(\delta\leq \delta_0\).
We start with the right-hand side inclusion of (\ref{IlluOvaloidEstimate}). Since \(|\conv[S, x^{\delta}]|_n-|S|_n=\delta|S|_n\), it follows that \(\langle x^{\delta}-x, u(x)\rangle=:\Delta\leq\Delta_0\). We conclude with Lemma \ref{CapHatEstimates}, 2., that 
\[
\delta|S|_n\geq \frac{(2\Delta)^{\frac{n+1}{2}}}{n(n+1)}\frac{|B_2^{n-1}|_{n-1}}{\kappa(x)^{1/2}}(1-\varepsilon)\quad,
\]
or, equivalently
\[
\Delta\leq \frac{n^{\frac{2}{n+1}}c(S,n)}{(1-\varepsilon)^{\frac{2}{n+1}}}\kappa(x)^{\frac{1}{n+1}}\delta^{\frac{2}{n+1}}\quad.
\]
Hence,
\begin{align}
\frac{\|x^{\delta}\|_2}{\|x\|_2}=\frac{\langle x^{\delta}, u(x)\rangle}{\langle x, u(x)\rangle}= 1+\frac{\Delta}{\langle x, u(x)\rangle}\leq 
1+\frac{n^{\frac{2}{n+1}}G(x)}{(1-\varepsilon)^{\frac{2}{n+1}}}\ \delta^{\frac{2}{n+1}}
\leq 1+\frac{n^{\frac{2}{n+1}}\GM}{(1-\varepsilon)^{\frac{2}{n+1}}}\ \delta^{\frac{2}{n+1}} , \notag
\end{align}
which proves the right-hand side of (\ref{IlluOvaloidEstimate}).
\par
\noindent
For the left-hand side inclusion,  let again \(x\in\partial S\). Then there is  \(x'\in\partial S\) such that \(x^{\delta}=x'+\Delta^{x'}(\delta)u(x')\), i.e., \(x'\) is the point of \(S\) with minimal  distance to \(x^{\delta}\). This point is unique since \(S\) is strictly convex. Note that \(\langle x-x', u(x')\rangle\leq 0\), i.e., \(\langle x, u(x')\rangle\leq \langle x', u(x')\rangle\) and therefore,
\[
\frac{\|x^{\delta}\|_2}{\|x\|_2}=\frac{\langle x^{\delta}, u(x')\rangle}{\langle x, u(x')\rangle}\geq \frac{\langle x'+\Delta^{x'}(\delta)u(x'), u(x')\rangle}{\langle x', u(x')\rangle}=1+\frac{\Delta^{x'}(\delta)}{\langle x', u(x')\rangle}\quad.
\]
Since \(\Delta^{x'}(\delta)\leq \Delta_0\),  Lemma \ref{CapHatEstimates} yields
\[
\frac{\Delta^{x'}(\delta)}{\langle x', u(x')\rangle}\geq \frac{n^{\frac{2}{n+1}}G(x')}{(1+\varepsilon)^{\frac{2}{n+1}}}\delta^{\frac{2}{n+1}}\geq \frac{n^{\frac{2}{n+1}}\Gm}{(1+\varepsilon)^{\frac{2}{n+1}}}\delta^{\frac{2}{n+1}}.
\] 
This establishes the  left-hand side inclusion  of (\ref{IlluOvaloidEstimate}). 
\hfill \(\Box\)
\vskip 3mm
\noindent
Note that  \(S^{\circ}\) is also a centrally symmetric convex body with \(C_+^2\)-boundary (see \cite{Hug}). An immediate corollary of Proposition \ref{UpperLowerEstimateIlluBody} is the following corollary.
\par
\noindent
\begin{cor}\label{DualIlluOvaloid} Suppose that \(S\) has \(C_+^2\)-boundary. Then there is a non-negative function \(\Psi\) with \(\lim_{\delta\rightarrow 0}\Psi(\delta)=0\) such that
\[
\left(1-(1+\Psi(\delta))\frac{\tilde{c}(S,n)}{\Gm}\delta^{\frac{2}{n+1}}\right)S\subseteq \ii[\delta]{S}\subseteq \left(1-(1-\Psi(\delta))\frac{\tilde{c}(S,n)}{\GM}\delta^{\frac{2}{n+1}}\right)S, 
\]
where \(\tilde{c}(S,n)=n^{\frac{2}{n+1}}c(S,n)c(S^{\circ},n)\).
\end{cor}
\textit{Proof.} The corollary follows immediately  from the following result which can be found in e.g., \cite{Hug}). For every \(x\in\partial S\) there exists a unique \(y\in\partial S^{\circ}\) with \(\langle x, y\rangle=1\) and in this case 
\[
\frac{\kappa_S(x)^{\frac{1}{n+1}}}{\langle x, u_S(x)\rangle}=\left[\frac{\kappa_{S^{\circ}}(y)^{\frac{1}{n+1}}}{\langle y, u_{S^{\circ}}(y)\rangle}\right]^{-1}\quad.
\]
\hfill \(\Box\)

\begin{thm}\label{MainTHMUpperBound}
\[
\limsup_{\delta\rightarrow 0}\frac{\d_S(\delta)-1}{\delta^{\frac{2}{n+1}}}\leq \GM-\Gm\quad.
\]
\end{thm}
\textit{Proof.} First assume that there is \(x\in\partial S\) with \(\kappa(x)=0\). Thus \(\GM-\Gm=\GM-0=\GM\). 
Note that
\[
\d_S(\delta)=\inf_{\delta'\geq 0}\d\left(S_{\delta},\ii{S}\right)
\leq \d\left(S_{\delta},  \ii[0]{S}\right)=\d\left(S_{\delta}, S\right)\quad.
\]
By Proposition \ref{LowerBoundFloatingBodyC2Case}, \((1-\GM\delta^{\frac{2}{n+1}}(1+o(1)))S\subseteq S_{\delta}\subseteq S\) and this yields
\[
\d\left(S_{\delta},S\right)\leq \frac{1}{1-\GM\delta^{\frac{2}{n+1}}(1+o(1))}=1+\GM\delta^{\frac{2}{n+1}}(1+o(1))\quad.
\]
Hence,
\[
\limsup_{\delta\rightarrow 0}\frac{\d_S(\delta)-1}{\delta^{\frac{2}{n+1}}}\leq \lim_{\delta\rightarrow 0}\GM(1+o(1))=\GM\quad.
\]
Assume now that \(S\) has \(C_+^2\)-boundary. Put \(\delta'=\frac{\GM\Gm}{\tilde{c}(S,n)}\delta^{\frac{2}{n+1}}\). By Proposition \ref{FloatingBodyOvaloid}, 
\[
\left(1-\GM\delta^{\frac{2}{n+1}}(1+o(1))\right)S\subseteq S_{\delta}\subseteq \left(1-\Gm\delta^{\frac{2}{n+1}}(1+o(1))\right)S . 
\]
By Corollary \ref{DualIlluOvaloid}, 
\[
\left(1-\GM\delta^{\frac{2}{n+1}}(1+o(1))\right)S\subseteq \ii{S}\subseteq \left(1-\Gm\delta^{\frac{2}{n+1}}(1+o(1))\right)S\quad.
\]
Therefore, a sufficient condition for \(\frac{1}{a}S_{\delta}\subseteq \ii{S}\) is 
\[
\frac{1}{a}\left(1-\Gm\delta^{\frac{2}{n+1}}(1+o(1))\right)\leq 1-\GM\delta^{\frac{2}{n+1}}(1+o(1))
\]
and a sufficient condition for \(\ii{S}\subseteq a S_{\delta}\) is
\[
1-\Gm\delta^{\frac{2}{n+1}}(1+o(1))\leq a\left(1-\GM\delta^{\frac{2}{n+1}}(1+o(1))\right)\quad.
\]
These two conditions are met, if one takes \(a=1+(\GM-\GM)\delta^{\frac{2}{n+1}}(1+o(1))\). Hence,
\[
\limsup_{\delta\rightarrow 0}\frac{\d_S(\delta)-1}{\delta^{\frac{2}{n+1}}}\leq \lim_{\delta\rightarrow 0}\frac{1+(\GM-\Gm)\delta^{\frac{2}{n+1}}(1+o(1))-1}{\delta^{\frac{2}{n+1}}}=\GM-\Gm
\]
\hfill \(\Box\)

\section{Lower Bounds}\label{SectionLowerBound}

We prove \(\liminf_{\delta\rightarrow 0}\frac{\d_S(\delta)-1}{\delta^{\frac{2}{n+1}}}\geq \GM-\Gm\) for \(S\) with \(C_+^2\)-boundary and for \(B_p^n\), \(2\leq p< \infty\). We also provide a lower bound for the case \(1<p<2\). Together with Theorem \ref{MainTHMUpperBound},  we get the following.
\begin{thm}\label{MainTHMLowerBound} Let \(S\subseteq \mathbb{R}^n\) be of class \(C_+^2\) or \(S=B_p^n\), \(2\leq p<\infty\). Then
\[
\lim_{\delta\rightarrow 0}\frac{\d_S(\delta)-1}{\delta^{\frac{2}{n+1}}}=\GM-\Gm\quad.
\]
\end{thm}
\vskip 3mm
\noindent
We need  another lemma.
\par
\noindent
\begin{lem}\label{FloatingIlluVectorsOvaloid} For every \(x\in\partial S\) the following holds.
\begin{enumerate}
\item There is a function \(\Phi_x\) with \(\lim_{\delta\rightarrow 0}\Phi_x(\delta)=0\) such that
\[
\left(1-(1+\Phi_x(\delta))G(x)\delta^{\frac{2}{n+1}}\right)x\in \partial S_{\delta}\quad.
\]
\item If \(S\) has \(C_+^2\)-boundary then there is a function \(\Psi_x\) with \(\lim_{\delta\rightarrow 0}\Psi_x(\delta)=0\) such that
\[
\left(1-(1+\Psi_x(\delta))\frac{\tilde{c}(S,n)}{G(x)}\delta^{\frac{2}{n+1}}\right)x\in\partial\left[\ii[\delta]{S}\right], 
\]
where \(\tilde{c}(S,n)\) is the constant defined in Lemma \ref{DualIlluOvaloid}.
\end{enumerate}
\end{lem} 
\vskip 3mm
\noindent
In order to prove this lemma,  we need two results. The first lemma is an immediate consequence of Lemmas 7 and 10 of \cite{SchuettWerner1990}.
\vskip 3mm
\noindent
\begin{lem}\label{FloatingBodyCurvature} Let \(K\subseteq\mathbb{R}^n\) be a convex body with \(0\) in its interior and \(x\in \partial K\) such that the Gauss curvature \(\kappa(x)\) exists. Then 
\[
\lim_{\delta\rightarrow 0}\frac{1}{\delta^{\frac{2}{n+1}}}\cdot \frac{\|x\|_2-\|x_{\delta}\|_2}{\|x\|_2}
=\frac{(n+1)^{\frac{2}{n+1}}}{2}\left(\frac{|K|_n}{|B_2^{n-1}|_{n-1}}\right)^{\frac{2}{n+1}}\frac{\kappa(x)^{\frac{1}{n+1}}}{\langle x, u(x)\rangle}=G(x), 
\]
where \(x_{\delta}\) is the unique point lying in the intersection of \(\partial K_{\delta}\) with the line segment \([0,x]\).
\end{lem}
\vskip 3mm
\noindent
The second lemma is a direct consequence of Lemma 3 of \cite{Werner94}.
\begin{lem}\label{IlluminationBodyCurvature}  Let \(K\subseteq\mathbb{R}^n\) be a convex body with \(0\) in its interior and \(x\in \partial K\) such that the Gauss curvature \(\kappa(x)\) exists. Then 
\[
\lim_{\delta\rightarrow 0}\frac{1}{\delta^{\frac{2}{n+1}}}\cdot \frac{\|x^{\delta}\|_2-\|x\|_2}{\|x\|_2}
=\frac{n^{\frac{2}{n+1}}(n+1)^{\frac{2}{n+1}}}{2}\left(\frac{|K|_n}{|B_2^{n-1}|_{n-1}}\right)^{\frac{2}{n+1}}\frac{\kappa(x)^{\frac{1}{n+1}}}{\langle x, u(x)\rangle}=n^{\frac{2}{n+1}}G(x)\quad,
\]
where \(x^{\delta}\) is the unique point on the boundary of \(K^{\delta}\) such that \(x\) lies on the line segment \([0,x^{\delta}]\).
\end{lem}

\textit{Proof of Lemma \ref{FloatingIlluVectorsOvaloid}.} It follows immediately from Lemma \ref{FloatingBodyCurvature} that for every \(x\in \partial S\) there is a function \(\Phi_x\) such that \(\lim_{\delta\rightarrow 0} \Phi_x(\delta)=0\) and
\[
\left(1-(1+\Phi_x(\delta))G(x)\delta^{\frac{2}{n+1}}\right)x\in \partial S_{\delta}\quad.
\]
\vskip 3mm
\noindent
The polar body \(S^{\circ}\) has \(C^2\)-boundary with everywhere positive Gauss curvature (see \cite{Hug}) and for every \(x\in \partial S\) there is a unique dual point \(y\), i.e. \(y\in \partial S^{\circ}\) such that  \(\langle x, y\rangle =1\) and 
\[
\frac{\kappa(x)^{\frac{1}{n+1}}}{\langle x, u(x)\rangle}
=\left[\frac{\kappa(y)^{\frac{1}{n+1}}}{\langle y, u(y)\rangle}\right]^{-1}\quad.
\]
By Lemma \ref{IlluminationBodyCurvature},
\[
\lim_{\delta\rightarrow 0} \frac{1}{\delta^{\frac{2}{n+1}}}\frac{\|y^{\delta}\|-\|y\|}{\|y\|}=n^{\frac{2}{n+1}}c(S^{\circ},n)\frac{\kappa(y)^{\frac{1}{n+1}}}{\langle y, u(y) \rangle}=n^{\frac{2}{n+1}}c(S^{\circ},n)\frac{\langle x, u(x)\rangle}{\kappa(x)^{\frac{1}{n+1}}}=\frac{\tilde{c}(S,n)}{G(x)}\quad.
\]
It follows that there is a function \(\tilde{\Psi}_y\) with \(\lim_{\delta\rightarrow 0} \tilde{\Psi}_y(\delta)=0\) such that
\[
\left(1+(1+\tilde{\Psi}_y(\delta))\frac{\tilde{c}(S,n)}{G(x)}\delta^{\frac{2}{n+1}}\right)y\in \left(S^{\circ}\right)^{\delta}\quad.
\]
Let \(x\in\partial S\) and \(y\in\partial S^{\circ}\) with \(\langle x,y\rangle=1\). We show that we may choose a function \(\Psi_x\)  such that  \(\lim_{\delta\rightarrow 0} \Psi_x(\delta)=0\) and such that 
\[
\left(1-(1+\Psi_x(\delta))\frac{\tilde{c}(S,n)}{G(x)}\delta^{\frac{2}{n+1}}\right)x\in \partial \left[\left(S^{\circ}\right)^{\delta}\right]^{\circ}=\ii[\delta]{S}\quad.
\]  
Let \(\lambda(\delta)\) be such that \(\lambda(\delta)x\in \ii[\delta]{S}\). Then
\[
1\geq\left\langle\lambda(\delta)x, \left(1+(1+\tilde{\Psi}_y(\delta))\frac{\tilde{c}(S,n)}{G(x)}\delta^{\frac{2}{n+1}}\right)y\right\rangle=\lambda(\delta) \left(1+(1+\tilde{\Psi}_y(\delta))\frac{\tilde{c}(S,n)}{G(x)}\delta^{\frac{2}{n+1}}\right)\quad.
\]
Hence we obtain that \(\lambda(\delta)\leq 1-(1+\Psi_x(\delta))\frac{\tilde{c}(S,n)}{G(x)}\delta^{\frac{2}{n+1}}\) for a suitable \(\Psi_x\). 
\par
\noindent
To establish  the opposite inequality, we use  techniques similar to the ones for the lower bound in the proof of [\cite{Werner94}, Lemma 3]. 
\par
\noindent
Translate and rotate \(S^{\circ}\) to a convex body \(K\) such that \(y\) is mapped to the origin and the outer normal is \(u_K(0)=-e_n\). Let \(f:\tau B_2^{n-1}\rightarrow \mathbb{R}\) be a parametrization of the boundary of \(K\) near the origin. Let \(\varepsilon>0\) and choose \(\eta>0\) such that
\[
\frac{1}{n}(1-\eta)^{n-1}-\frac{1}{n+1}(1+\eta)^{n-1}\geq \frac{1}{n(n+1)}\frac{1}{(1+\varepsilon)^{\frac{n+1}{2}}}\quad.
\] 
Then there exists \(\Delta_0>0\) such that for every \(0\leq\Delta\leq\Delta_0\) 
\[
(1-\eta)\sqrt{2\Delta} \mathcal{E}\subseteq \{z\in\tau B_2^{n-1}: f(z)\leq \Delta\}\subseteq (1+\eta)\sqrt{2\Delta}\mathcal{E}
\]
where \(\mathcal{E}=\{z\in\mathbb{R}^{n-1}:\langle Hf(0)z,z \rangle\leq 1\}\) is the indicatrix of Dupin (see, e.g.,  \cite{ReisnerSchuettWerner}, \cite{Schuett:2003}).
We conclude that for every \(\zeta\in\mathbb{R}^{n-1}\)
\[
\mathrm{conv}[K,\zeta-\Delta e_n]\backslash K\supseteq
\mathrm{conv}[\Delta e_n+(1-\eta)\sqrt{2\Delta}\mathcal{E}, \zeta-\Delta e_n]\backslash K\quad. 
\]
We compute
\begin{eqnarray*}
&& \hskip -5mm \left|\mathrm{conv}[K,\zeta-\Delta e_n]\backslash K\right|_n\geq
\left|\mathrm{conv}[\Delta e_n+(1-\eta)\sqrt{2\Delta}\mathcal{E}, \zeta-\Delta e_n]\backslash K\right|_n \\
&& \hskip 5mm\geq|\mathrm{conv}[\Delta e_n+(1-\eta)\sqrt{2\Delta}\mathcal{E}, \zeta-\Delta e_n]|_n-|\{v\in K: v_n \leq \Delta\}|_n \quad.
\end{eqnarray*}
We have
\[
|\mathrm{conv}[\Delta e_n+(1-\eta)\sqrt{2\Delta}\mathcal{E}, \zeta-\Delta e_n]|_n=(1-\eta)^{n-1}\frac{(2\Delta)^{\frac{n+1}{2}}}{n}|\mathcal{E}|_{n-1}
\]
and 
\begin{eqnarray*}
|\{v\in K: v_n \leq \Delta\}|_n &\leq& \int_0^{\Delta}|\{z\in\mathbb{R}^{n-1}:f(z)\leq t\}|_{n-1}\mathrm{d}
t\\
&\leq &(1+\eta)^{n-1}\int_0^{\Delta} (2t)^{\frac{n-1}{2}}|\mathcal{E}|_{n-1}\mathrm{d}t
=(1+\eta)^{n-1}\frac{(2\Delta)^{\frac{n+1}{2}}}{n+1}|\mathcal{E}|_{n-1}\quad.\notag
\end{eqnarray*}
Therefore,
\begin{align}
&\left|\mathrm{conv}[K,\zeta-\Delta e_n]\backslash K\right|_n\geq (2\Delta)^{\frac{n+1}{2}}|\mathcal{E}|_{n-1}\left(\frac{(1-\eta)^{n-1}}{n}-\frac{(1+\eta)^{n-1}}{n+1}\right)
\geq \frac{(2\Delta)^{\frac{n+1}{2}}}{n(n+1)}\frac{|\mathcal{E}|_{n-1}}{(1+\varepsilon)^{\frac{n+1}{2}}}\quad.\notag
\end{align}
Suppose that \(\delta_0>0\) is sufficiently small such that for every \(0\leq \delta\leq \delta_0\) we have
\[
\frac{(2\Delta_0)^{\frac{n+1}{2}}}{n(n+1)}\frac{|\mathcal{E}|_{n-1}}{(1+\varepsilon)^{\frac{n+1}{2}}}\geq \delta|K|_n\quad.
\]
Hence, if \(0\leq \delta \leq \delta_0\)  and \(\Delta\geq 0\) is chosen such that
\[
\left|\mathrm{conv}[K,\zeta-\Delta e_n]\backslash K\right|_n=\delta|K|_n
\]
it follows that \(\Delta\leq \Delta_0\) and we conclude that
\[
\delta|K|\geq \frac{(2\Delta)^{\frac{n+1}{2}}}{n(n+1)}\frac{|\mathcal{E}|_{n-1}}{(1+\varepsilon)^{\frac{n+1}{2}}}\quad,
\]
or, equivalently
\[
\Delta\leq (1+\varepsilon)\frac{1}{2}\left[n(n+1)\frac{|K|_n}{|B_2^{n-1}|_{n-1}}\right]^{\frac{2}{n+1}}\kappa_K(0)^{\frac{1}{n+1}}\delta^{\frac{2}{n+1}}\quad.
\]
Hence, 
\[
K^{\delta}\subseteq \ll\xi\in\mathbb{R}^n: \xi_n\geq -(1+\varepsilon)\frac{1}{2}\left[n(n+1)\frac{|K|_n}{|B_2^{n-1}|_{n-1}}\right]^{\frac{2}{n+1}}\kappa_K(0)^{\frac{1}{n+1}}\delta^{\frac{2}{n+1}}\rr\quad.
\]
Thus,
\[
\left(S^{\circ}\right)^{\delta}\subseteq \{v\in\mathbb{R}^n: \langle v-y, u(y) \rangle\leq (1+\varepsilon)n^{\frac{2}{n+1}}c(S^{\circ},n)\kappa(y)^{\frac{1}{n+1}}\delta^{\frac{2}{n+1}}\}\quad.
\]
Let \(x\in \partial S\) be the unique point such that \(\langle x,y \rangle=1\). It follows that for \(0\leq \delta \leq \delta_0\) 
\begin{align}
\lambda(\delta)\geq & \left(1+(1+\varepsilon)n^{\frac{2}{n+1}}c(S^{\circ},n)\frac{\kappa(y)^{\frac{1}{n+1}}}{\langle y, u(y)\rangle } \delta^{\frac{2}{n+1}}\right)^{-1}
\geq 1-(1+\varepsilon)n^{\frac{2}{n+1}}c(S^{\circ},n)\frac{\kappa(y)^{\frac{1}{n+1}}}{\langle y, u(y)\rangle } \delta^{\frac{2}{n+1}}\notag\\
=&1-(1+\varepsilon)n^{\frac{2}{n+1}}c(S^{\circ},n)\frac{\langle x, u(x)\rangle}{\kappa(x)^{\frac{1}{n+1}}}\delta^{\frac{2}{n+1}}=1-(1+\varepsilon)\frac{\tilde{c}(S,n)}{G(x)}\delta^{\frac{2}{n+1}}. 
\end{align}
This proves that we may choose a function \(\Psi_x\) such that \(\lim_{\delta\rightarrow 0}\Psi_x(\delta)=0\) and 
\[
\left(1-(1+\Psi_x(\delta))\frac{\tilde{c}(S,n)}{G(x)}\right)x\in \partial \left[\ii[\delta]{S}\right]\quad.
\]
\hfill \(\Box\)

\subsection{The $C_+^2$-Case}
\vskip 3mm
\noindent
\begin{prop}\label{LowerBoundOvaloid} Suppose that \(S\) has \(C_+^2\)-boundary. Then
\[
\liminf_{\delta\rightarrow 0}\frac{\d_S(\delta)-1}{\delta^{\frac{2}{n+1}}}\geq \GM-\Gm\quad.
\]
\end{prop}
\noindent
\textit{Proof.} For \(x\in\partial S\) apply Lemma \ref{FloatingIlluVectorsOvaloid} and put  \(\delta''=\tilde{c}(S,n)^{\frac{n+1}{2}}\delta'\), i.e., 
\[
(1-(1+\Omega_x(\delta''))G(x)^{-1}\delta''^{\frac{2}{n+1}})x\in\partial\ii{S}, 
\]
where \(\Omega_x(\delta'')=\Psi_x(\delta')\). Let \(x_0\in\partial S\) be such that \(\Gm=\min_{x\in\partial S}G(x)=G(x_0)\) and let \(x_1\in \partial S\) be such that  \(\GM=\max_{x\in\partial S}G(x)=G(x_1)\). Let \(a\geq 1\) be such that \(\frac{1}{a}S_{\delta}\subseteq \left[(S^{\circ})^{\delta'}\right]^{\circ}\subseteq a S_{\delta}\). 
We conclude that 
\begin{align}
a\geq \frac{1-(1+\Phi_{x_0}(\delta))\Gm\delta^{\frac{2}{n+1}}}{1-(1+\Omega_{x_0}(\delta''))\Gm^{-1}{\delta''}^{\frac{2}{n+1}}}\label{Estimate1}
\end{align}
and 
\begin{align}
a\geq \frac{1-(1+\Omega_{x_1}(\delta''))\GM^{-1}{\delta''}^{\frac{2}{n+1}}}{1-(1+\Phi_{x_1}(\delta))\GM\delta^{\frac{2}{n+1}}}\quad.\label{Estimate2}
\end{align}
The assumption that \(\delta''\geq (\Gm\GM)^{\frac{n+1}{2}}\delta\) and inequality (\ref{Estimate1}) lead to 
$$
\liminf_{\delta>0}\frac{a-1}{\delta^{\frac{2}{n+1}}}\geq \GM-\Gm. 
$$
The assumption \(\delta''\leq (\Gm\GM)^{\frac{n+1}{2}}\delta\) and (\ref{Estimate2}) lead to 
$$
\liminf_{\delta>0}\frac{a-1}{\delta^{\frac{2}{n+1}}}\geq \GM-\Gm
$$ 
as well. We conclude
that
\[
\liminf_{\delta\rightarrow 0}\frac{\d_S(\delta)-1}{\delta^{\frac{2}{n+1}}}\geq \GM-\Gm\quad .
\]
\hfill \(\Box\)

\subsection{The $B_p^n$-Case}
In this section we prove the following.
\begin{thm}\label{MainTheoremLP} For every \(n\geq 2\) and \(2<p<\infty\) we have
\[
\lim_{\delta\rightarrow 0}\frac{\d_p(\delta)-1}{\delta^{\frac{2}{n+1}}}=G((n^{-1/p},\dots,n^{-1/p}))=\GM=\GM-\Gm\quad.
\]
\end{thm}
\vskip 3mm
\noindent
\begin{rem} By Lemma \ref{lpCurvatureNormal} below, the value of \(\GM\) is given by
\[
G((n^{-1/p},\dots,n^{-1/p}))=c(B_p^n,n)(p-1)^{\frac{n-1}{n+1}}n^{\frac{2n+p}{(n+1)p}}\quad.
\]
\end{rem}
\vskip 3mm
\noindent
Furthermore, we give a lower bound for the case \(1<p<2\). First, we need some preparatory definitions and lemmas.
\begin{defn} Let \(1\leq p\leq \infty\). The cap \(C_p^n(\Delta)\) of \(B_p^n\) of height \(0\leq \Delta\leq 1\) is defined by
\[
C_p^n(\Delta)=\{ x\in B_p^n: x_1\geq 1-\Delta\}.
\]
The hat \(H_p^n(\Delta)\) of height \(\Delta\geq 0\) is defined by
\[
H_p^n(\Delta)=\mathrm{conv}[B_p^n,(1+\Delta)e_1]\backslash B_p^n.
\]
\end{defn}

\begin{lem}\label{lrKalottenvolumen}
Let \(1\leq p <\infty\) and \(0\leq \Delta\leq 1\). Then
\[
|C_p^n(\Delta)|_n=\frac{1}{n-1+p}|B_p^{n-1}|_{n-1}(p\Delta)^{\frac{n-1+p}{p}}(1-\phi_p^n(\Delta))\quad,
\]
where \(\phi_p^n\) is a function with \(\lim\limits_{\Delta\rightarrow 0}\phi_p^n(\Delta)=0\).
\end{lem}
\textit{Proof.} By Taylor's theorem,  for every \(0<s<1\) there is  \(0<\sigma(s)<s\) such that \((1-s)^p=1-p(1-\sigma(s))^{p-1}s\). This yields
\begin{align}
\int_0^{\Delta} \left(1-(1-s)^p\right)^{\frac{n-1}{p}}\mathrm{d}s=\int_0^{\Delta} \left(p(1-\sigma(s))^{p-1}s\right)^{\frac{n-1}{p}}\mathrm{d}s
\leq  \int_0^{\Delta}(ps)^{\frac{n-1}{p}}\mathrm{d}s= \frac{1}{n-1+p}(p\Delta)^{\frac{n-1+p}{p}}\notag
\end{align}
and 
\begin{align}
\int_0^{\Delta} \left(p(1-\sigma(s))^{p-1}s\right)^{\frac{n-1}{p}}\mathrm{d}s
\geq (1-\Delta)^{\frac{p-1}{p}(n-1)}\int_0^{\Delta}(ps)^{\frac{n-1}{p}}\mathrm{d}s
=(1-\Delta)^{\frac{p-1}{p}(n-1)} \frac{1}{n-1+p}(p\Delta)^{\frac{n-1+p}{p}}\quad.\notag
\end{align}
Hence, we obtain by Cavalieri's principle
\[
|C_p^n(\Delta)|_n=\frac{1}{n-1+p}|B_p^{n-1}|_{n-1}(p\Delta)^{\frac{n-1+p}{p}}(1-\phi_p^n(\Delta)), 
\]
with \(0\leq \phi_p^n(\Delta)\leq 1-(1-\Delta)^{\frac{p-1}{p}(n-1)}\). 

\hfill \(\Box\)
\vskip 3mm
\noindent
\begin{lem}\label{lrHutvolumen}
Let \(1<p<\infty\) and \(0\leq \Delta\leq 1\). Then
\[
|H_p^n(\Delta)|_n=\frac{p-1}{n(n-1+p)}|B_p^{n-1}|_{n-1}\left(\frac{p}{p-1}\Delta\right)^{\frac{n-1+p}{p}}(1+\psi_p^n(\Delta))\quad,
\]
where \(\psi_p^n\) is a function with \(\lim\limits_{\Delta\rightarrow 0}\psi_p^n(\Delta)=0\).
\end{lem}
\textit{Proof.} Put \(f(t)=(1-t^p)^{\frac{1}{p}}\), \(-1<t<1\). The tangential function \(T_{t_0}(t)\) of \(f\) at \(t_0\) is given by
\[
T_{t_0}(t)=f(t_0)+f'(t_0)(t-t_0)=\left(1-t_0^p\right)^{\frac{1}{p}-1}(1-t_0^{p-1}t)\quad.
\]
If \(T_{t_0}(1+\Delta)=0\), then \(t_0=\frac{1}{(1+\Delta)^{\frac{1}{p-1}}}\). It follows that the boundary of the cap of \(B_p^n\) touches \(B_p^n\) at height \(\frac{1}{(1+\Delta)^{\frac{1}{p-1}}}\) with respect to the direction \(e_1\). The hat volume is given by the difference of the cone volume and the corresponding cap volume, 
\begin{align}
|H_p^n(\Delta)|_n=\frac{1}{n}\left(1+\Delta-\frac{1}{(1+\Delta)^{\frac{1}{p-1}}}\right)\left(1-\frac{1}{(1+\Delta)^{\frac{p}{p-1}}}\right)^{\frac{n-1}{p}}|B_p^{n-1}|_{n-1}
- \left| C_p^n\left(1-\frac{1}{(1+\Delta)^{\frac{1}{p-1}}}\right)\right|_n. \notag
\end{align}
Hence, 
\begin{align}
\lim\limits_{\Delta\rightarrow 0} \frac{1}{\Delta^{\frac{n-1+p}{p}}}|H_p^n(\Delta)|_n
=&\frac{1}{n}\left(\frac{p}{p-1}\right)^{\frac{n-1+p}{p}}|B_p^{n-1}|_{n-1}\notag
-\frac{1}{n-1+p}\left(\frac{p}{p-1}\right)^{\frac{n-1+p}{p}}|B_p^{n-1}|_{n-1}\notag\\
=&\frac{p-1}{n(n-1+p)}\left(\frac{p}{p-1}\right)^{\frac{n-1+p}{p}}|B_p^{n-1}|_{n-1}\quad . \notag
\end{align}
\hfill \(\Box\)

\vspace{1cm}
\noindent
The following formulas for the Gauss curvature and the normal directions of \(B_p^n\) can be found in e.g. \cite{KoldobskyRyaboginZvavitch, SchuettWerner1990}.
\begin{lem}\label{lpCurvatureNormal} The Gauss curvature at \(x\in\partial B_p^n\) is given by
\[
\kappa(x)=(p-1)^{n-1}\frac{\prod_{i=1}^n |x_i|^{p-2}}{\left(\sum_{i=1}^n |x_i|^{2p-2}\right)^{\frac{n+1}{2}}}
\]
and the normal is given by
\[
u(x)=\frac{1}{(\sum_{i=1}^n |x_i|^{2p-2})^{\frac{1}{2}}}\left(\mathrm{sign} (x_i)|x_i|^{p-1}\right)_{i=1}^n\quad,
\]
if \(2\leq p <\infty\) or if \(1<p<2\) and all components of \(x\) are not equal to zero.
\end{lem}  

\vskip 3mm
\noindent
\begin{lem}\label{FloatingIlluminationVectors} Let \(1<p<\infty\). 
\begin{enumerate}
\item There are functions \(\Phi_1\), \(\Phi_2\), depending only on \(n\) and \(p\), such that \(\lim\limits_{\delta\rightarrow 0}\Phi_1(\delta)=0\), \(\lim\limits_{\delta\rightarrow 0}\Phi_2(\delta)=0\) and
\begin{align}
\pm(1-(1+\Phi_1(\delta))c_1\delta^{\frac{p}{n-1+p}})e_{i}
\in & \partial(B_p^n)_{\delta},\quad 1\leq i \leq n,
\notag\\
\left(1-(1+\Phi_2(\delta))c_2\delta^{\frac{2}{n+1}}\right)
\frac{1}{\sqrt[p]{n}}\sum_{i=1}^n \pm e_i\in & \partial (B_p^n)_{\delta}\notag
\end{align}
where
\begin{align*}
c_1=\frac{(n-1+p)^{\frac{p}{n-1+p}}}{p}\left(\frac{|B_p^n|_n}{|B_p^{n-1}|_{n-1}}\right)^{\frac{p}{n-1+p}},& \qquad c_2=G_{B_p^n}((n^{-1/p},\dots,n^{-1/p}))\quad.
\end{align*}
\item There are functions \(\Psi_1\), \(\Psi_2\), depending only on \(n\) and \(p\), such that \(\lim\limits_{\delta\rightarrow 0}\Psi_1(\delta)=0\), \(\lim\limits_{\delta\rightarrow 0}\Psi_2(\delta)=0\) and
\begin{align}
\pm(1-(1+\Psi_1(\delta))c_3\delta^{\frac{p}
{n-1+p}})e_{i}\in & \partial \ii[\delta]{B_{p'}^n},\quad 1\leq i\leq n,\notag\\
\left(1-(1+\Psi_2(\delta))c_4\delta^{\frac{2}{n+1}}
\right)\frac{1}{\sqrt[p]{n}}\sum_{i=1}^n \pm e_i\in & 
\partial\ii[\delta]{B_{p'}^n}\quad,
\notag
\end{align}
where \(\frac{1}{p}+\frac{1}{p'}=1\) and
\begin{align*}
c_3=\left[\frac{n(n-1+p')}{p'-1}\cdot\frac{|B_{p'}^n|_n}{|B_{p'}^{n-1}|_{n-1}}\right]^{\frac{p'}{n-1+p'}}\cdot\frac{p'-1}{p'},& \qquad c_4=n^{\frac{2}{n+1}}G_{B_{p'}^n}((n^{-1/p'},\dots,n^{-1/p'}))\quad.
\end{align*}
\end{enumerate}
\end{lem}
\textit{Proof.} Part 1. is an immediate consequence of Lemmas \ref{lrKalottenvolumen} and \ref{FloatingBodyCurvature}. For Part 2.,  note that by Lemma \ref{lrKalottenvolumen}
\[
\pm(1+(1+\psi_1(\delta))c_3\delta^{\frac{p}{n-1+p}})e_i\in \partial\left[(B_p^n)^{\delta}\right].
\]
By Lemma \ref{IlluminationBodyCurvature}
\[
\left(1+(1+\psi_2(\delta))c_4\delta^{\frac{2}{n+1}}
\right)\frac{1}{\sqrt[p]{n}}\sum_{i=1}^n \pm e_i\in  
\partial \left[(B_p^n)^{\delta}\right]\quad,
\] 
for some functions \(\psi_1\), \(\psi_2\) with \(\lim_{\delta\rightarrow 0}\psi_1(\delta)=0\) and \(\lim_{\delta\rightarrow 0}\psi_2(\delta)=0\). Furthermore, note that by symmetry \(\pm e_i\) is an outer normal of \((B_p^n)^{\delta}\) at
\(
\pm(1+(1+\psi_1(\delta))c_3\delta^{\frac{p}{n-1+p}})e_i
\) 
and \((n^{-1/2},\dots,n^{-1/2})\) is an outer normal of \((B_p^n)^{\delta}\) at
\(
\left(1+(1+\psi_2(\delta))c_4\delta^{\frac{2}{n+1}}
\right)\frac{1}{\sqrt[p]{n}}\sum_{i=1}^n \pm e_i
\).
Similar to the proof of Part 2. of Lemma \ref{FloatingIlluVectorsOvaloid},  one can show that
\[
\pm(1+(1+\psi_1(\delta))c_3\delta^{\frac{p}{n-1+p}})^{-1}e_i\in\ii[\delta]{B_{p'}^n}
\]
and
\[
\left(1+(1+\psi_2(\delta))c_4\delta^{\frac{2}{n+1}}
\right)^{-1}\frac{1}{\sqrt[p']{n}}\sum_{i=1}^n \pm e_i\in\ii[\delta]{B_{p'}^n}\quad.
\]
From this, we conclude.
\hfill \(\Box\)

\vskip 3mm
\noindent

\begin{prop}
Let \(2\leq p<\infty\) and let \(c_2\) be the constant of Lemma \ref{FloatingIlluminationVectors}. Then we have 
\[
\liminf_{\delta\rightarrow 0}\frac{\d_p(\delta)-1}{\delta^{\frac{2}{n+1}}}\geq c_2=\GM=G((n^{-1/p},\dots,n^{-1/p}))\quad.
\]
\end{prop}
\textit{Proof.} Fix \(\delta'\geq 0\) and let \(a_{\delta}\geq 1\) be such that
\[
\frac{1}{a_{\delta}} (B_p^n)_{\delta}\subseteq \ii{B_p^n}\subseteq a_{\delta} (B_p^n)_{\delta}\quad.
\]
From Lemma \ref{FloatingIlluminationVectors} we deduce that
\[
a_{\delta}\geq \max \left[\frac{1-c_4\delta'^{\frac{2}{n+1}}(1+\Psi_4(\delta'))}{1-c_2\delta^{\frac{2}{n+1}}(1+\Phi_2(\delta))}, \frac{1-c_1\delta^{\frac{p}{n-1+p}}(1+\Phi_1(\delta))}{1-c_3\delta'^{\frac{p'}{n-1+p'}}(1+\Psi_3(\delta'))}\right]
\]
Since 
\(\frac{p'}{n-1+p'}<\frac{2}{n+1}<\frac{p}{n-1+p}\), we may choose a fixed constant \(1<\alpha<\frac{2}{n+1}\cdot\frac{n-1+p'}{p'}\). As \(\alpha \frac{p'}{n-1+p'}<\frac{2}{n+1}\),  
we get for \(\delta'\geq  \delta^{\alpha}\)  that
\begin{align}
a_{\delta}\geq &\frac{1-c_1\delta^{\frac{p}{n-1+p}}(1+s_1(\delta))}{1-c_3\delta^{\alpha \frac{p'}{n-1+p'}}(1+s_3(\delta'))}=1+c_3\delta^{\alpha \frac{p'}{n-1+p'}}(1+o_{\delta}(1)+o_{\delta'}(1))\notag\\
\geq & 1+c_2\delta^{\frac{2}{n+1}}(1+o_{\delta}(1)+o_{\delta'}(1))\quad.\notag
\end{align}
If \(\delta'\leq \delta^{\alpha}\) it follows that
\[
a_{\delta}\geq \frac{1-c_4\delta^{\alpha \frac{2}{n+1}}(1+s_4(\delta'))}{1-c_2\delta^{\frac{2}{n+1}}(1+s_2(\delta))}=1+c_2\delta^{\frac{2}{n+1}}(1+o_{\delta}(1)+o_{\delta'}(1))\quad.
\]
We may assume without loss of generality that \(\delta'\) tends to zero as \(\delta\) tends to zero and hence,
\[
\lim_{\delta\rightarrow 0}\frac{a_{\delta}-1}{\delta^{\frac{2}{n+1}}}\geq c_2\quad.
\]
The last step is to prove that \(\GM=G((n^{-1/p},\dots,n^{-1/p}))\).
By Lemma \ref{lpCurvatureNormal}, 
\[
\frac{\kappa(x)^{\frac{1}{n+1}}}{\langle x, u(x)\rangle}=(p-1)^{\frac{n-1}{n+1}}\prod_{i=1}^n|x_i|^{\frac{p-2}{n+1}}\quad.
\]
The arithmetic geometric means inequality yields
\[
\prod_{i=1}^n |x_i|^{\frac{p-2}{n+1}}=\left(\sqrt[n]{\prod_{i=1}^n |x_i|^p}\right)^{\frac{n}{n+1}\cdot\frac{p-2}{p}}\leq \left(\frac{1}{n}\left(\sum_{i=1}^n |x|_i^p\right)\right)^{\frac{n}{n+1}\cdot\frac{p-2}{p}}
=n^{-\frac{n}{n+1}\cdot\frac{p-2}{p}}
\] , 
with equality if and only if all \(|x_i|\) are equal to \(n^{-1/p}\).
\hfill \(\Box\)
\vskip 3mm
\noindent
\begin{prop} Let \(1<p<2\) and let \(c_1\) be the constant from Lemma \ref{FloatingIlluminationVectors}. We have 
\[
\liminf_{\delta\rightarrow 0}\frac{\d_p(\delta)-1}{\delta^{\frac{p}{n-1+p}}}\geq c_1=\frac{(n-1+p)^{\frac{p}{n-1+p}}}{p}\left(\frac{|B_p^n|_n}{|B_p^{n-1}|_{n-1}}\right)^{\frac{p}{n-1+p}}\quad.
\] 
\end{prop}
\textit{Proof.} Similarly as in the previous proposition one argues that
\[
a_{\delta}\geq \max\left[\frac{1-c_2\delta^{\frac{2}{n+1}}(1+o_{\delta}(1))}{1-c_4\delta'^{\frac{2}{n+1}}(1+o_{\delta'}(1))}, \frac{1-c_3\delta'^{\frac{p'}{n-1+p'}}}{1-c_1\delta^{\frac{p}{n-1+p}}}\right]\quad.
\]
Note that we have \(\frac{p}{n-1+p}< \frac{2}{n+1}<\frac{p'}{n-1+p'}\). We fix an \(\alpha>0\) with
\[
\frac{p}{n-1+p}\cdot\frac{n-1+p'}{p'}<\alpha < \frac{p}{n-1+p}\cdot \frac{n+1}{2}
\]
and we conclude as in the previous proof by considering the two cases \(\delta'\leq\delta^{\alpha}\) and \(\delta'\geq \delta^{\alpha}\). 

\hfill \(\Box\)

\subsection*{Acknowledgement}
The first author likes to thank the Department of Mathematics at Case Western Reserve University, Cleveland, for their hospitality during his research stay in 2015/2016.
Both authors want to thank the Mathematical Science Research Institute, Berkley. It was during a stay there when the paper was completed.

\vskip 4mm
\noindent
Olaf Mordhorst\\
{\small Institut f\"ur diskrete Mathematik und Geometrie}\\
{\small Technische Universit\"at Wien}\\
{\small 1040 Wien, Austria}\\
{\small \tt olaf.mordhorst@tuwien.ac.at}\\ \\

\vskip 3mm
\noindent
Elisabeth M. Werner\\
{\small Department of Mathematics \ \ \ \ \ \ \ \ \ \ \ \ \ \ \ \ \ \ \ Universit\'{e} de Lille 1}\\
{\small Case Western Reserve University \ \ \ \ \ \ \ \ \ \ \ \ \ UFR de Math\'{e}matique }\\
{\small Cleveland, Ohio 44106, U. S. A. \ \ \ \ \ \ \ \ \ \ \ \ \ \ \ 59655 Villeneuve d'Ascq, France}\\
{\small \tt elisabeth.werner@case.edu}\\ \\

\end{document}